\documentclass{ws-ijm}

\begin{document}

\markboth{S. Zhang and Y.-Z. Zhang}
{Pointed Hopf algebras with classical  Weyl groups}

\catchline{}{}{}{}{}

\title{POINTED HOPF ALGEBRAS WITH CLASSICAL WEYL GROUPS}

\author{SHOUCHUAN ZHANG}

\address{Department  of Mathematics, Hunan University,
 Changsha  410082,  P.R. China }

\author{ YAO-ZHONG ZHANG}

\address{School of Mathematics and Physics, The University of Queensland,
 Brisbane 4072, Australia}

\maketitle

\begin{abstract}
We prove that  Nichols algebras of irreducible
Yetter-Drinfeld modules over classical Weyl groups $A \rtimes \mathbb
S_n$  supported by $\mathbb S_n$ are infinite dimensional, except in three cases. We give
necessary and sufficient conditions for Nichols algebras of
Yetter-Drinfeld modules over classical Weyl groups $A \rtimes \mathbb S_n$
supported by  $A$ to be finite dimensional.
\end{abstract}



\newtheorem{Proposition}{Proposition}[section]
\newtheorem{Theorem}[Proposition]{Theorem}
\newtheorem{Definition}[Proposition]{Definition}
\newtheorem{Corollary}[Proposition]{Corollary}
\newtheorem{Lemma}[Proposition]{Lemma}
\newtheorem{Example}[Proposition]{Example}
\newtheorem{Remark}[Proposition]{Remark}

\maketitle 

\numberwithin{equation}{section}

\section{Introduction}\label {s0}

This paper is a contribution to the  classification of finite
dimensional complex pointed Hopf algebras with non abelian group $G$.
 It is known that the first necessary step is to consider the
 Nichols algebras associated to the irreducible Yetter-Drinfeld
 modules over $G$ and decide if they are finite dimensional or not \cite{AS98,AS02}. 
These irreducible Yetter-Drinfeld
 modules are easy to classify: they are in one-to -one correspondence
 with pairs $({\mathcal C}, \rho)$ where ${\mathcal C}$ is a
 conjugacy class of $G$ and $\rho$ is an irreducible representation
 of the centralizer of an element in $\mathcal C$. In short, it is
 necessary  to see if the dimension of the corresponding Nichols
 algebra $\mathfrak B({\mathcal C}, \rho)$ is finite or not. As
 pointed out early by Gra\~na \cite{Gr00}, one may start attacking this problem by looking at Nichols
 subalgebras of $\mathfrak B({\mathcal C}, \rho)$.

 This paper deals specifically with the case when  $G$ is the Weyl
 group of a simple Lie algebra. The most prominent example is the
 symmetric group  $\mathbb S_n.$ For this group, a fairly complete
 analysis is presented in \cite {AFGV}, as a culmination of the
  series of papers of  \cite {AZ07,AFZ08}. The outcome is that all $\mathfrak B({\mathcal C},
  \rho)$ have  infinite dimension, except for a small list of examples
  when $n\le 6$ and remarkable cases corresponding to ${\mathcal C}$
  = the class of transpositions, and $\rho =$ the restriction of representation sign
  or a closely related one, whose dimension is still unknown. It is natural to guess that the conjugacy class of involutions in other Weyl
  groups would also  be distinguished.

  From Chapters IV, V, VI of \cite {Bo68}, $(\mathbb Z_2)^n \rtimes \mathbb S_n$
is isomorphic to the Weyl groups $W(B_n)$ and $W(C_n)$ of $B_n$ and
$C_n$, where $n>2$. If $A= \{a \in
     (\mathbb Z_2)^n \mid \ a = ({a_1}, {a_2}, \cdots, {a_n}) $  $\hbox {
with   } a_1 +a_2 + \cdots + a_n \equiv 0  \ (\ {\rm mod } \ 2) \}$,
then $A \rtimes \mathbb S_n$ is isomorphic to the Weyl group
$W(D_n)$ of $D_n$, where $n>3$. When $A= \{a \in
     (\mathbb Z_2)^n \mid \ a = ({a_1}, {a_2}, \cdots, {a_n}) $  $\hbox {
with  all } a_i =0 \}$, $A \rtimes \mathbb S_n$ are isomorphic to the Weyl group of $A_{n-1}$, where $n>1$.

 Note that $\mathbb{S}_n$ acts on $A$ as follows. For any $a\in A$ with $a = ({a_1}, {a_2},
\cdots, {a_n})$ and $h\in \mathbb{S}_n$, 
 \begin {eqnarray} \label {action}
h \cdot a =:( {a_{h^{-1}(1)}}, {a_{h^{-1}(2)}}, \cdots,{a_{h^{-1}(n)}} ).
\end {eqnarray}
It is clear that 
\begin {equation}\label {e2.10} (a, \sigma) ^{-1} = (-(a_{\sigma (1)}, a_{\sigma (2)}, \cdots, a_{\sigma (n)}),  
      \sigma ^{-1}) = (-\sigma ^{-1}(a),  \sigma ^{-1}) , 
\end  {equation}
\begin {equation}\label {e2.11}
(b, \tau) (a, \sigma) (b, \tau)^{-1} = (b + \tau (a) - \tau \sigma \tau ^{-1} (b), \tau \sigma \tau ^{-1}).
\end {equation}

If $H$ is a subgroup of a group $G$ and $\mathcal O$ is a conjugacy class in $G$, 
then we say that $\mathcal O$ is supported by $H$ if ${\mathcal O} \cap H \not=  \emptyset$.

Without specification,   $ G=:A\rtimes  \mathbb{S}_n$  with  $ A=: (\mathbb Z_2)^n $ or $A \subseteq (\mathbb Z_2)^n $. 
Let $\chi _2$ denote the character of ${\mathbb Z}_2$ with order $2$, i.e. $\chi _2 (1) =-1.$
Let $\sigma\in \mathbb S_n$ be of type $(1^{\lambda_1}, 2^{\lambda_2}, \dots,
n^{\lambda_n})$. If   $\rho \in \widehat {G^\sigma }$, then $\rho $
may be represented as  $\rho = \theta_{\chi , \mu } = (\chi \otimes \mu )\uparrow ^{G^\sigma} _{(G^\sigma
)_\chi} \in \widehat {A^\sigma \rtimes \mathbb S_n ^\sigma }$ with
$\chi = (\chi _2 ^{b_1} \otimes \cdots \chi_2 ^{b_n})  \in \widehat
{A^{\sigma}}$ and $\mu \in \widehat {({\mathbb S}_n^\sigma )_\chi}$, where
$\mu =\otimes _{ 1\le i \le n } \mu _i $ with $\mu _i \in
\widehat { \mathbb S_{Y_i}^{\sigma_i}}$.  Here and below $W_\chi$ and  $Y_i$  
are defined  in subsections \ref {s1.1} and \ref {s1.3}, respectively.
Assume that $M = M({\mathcal O}_{\sigma _1}, \rho ^{(1)})\oplus M({\mathcal
O}_{\sigma _2}, \rho ^{(2)}) \oplus \cdots \oplus M({\mathcal
O}_{\sigma_m}, \rho ^{(m)})$ is a reducible {\rm YD} module over ${\mathbb C } G$. 

The main results in this paper are summarized in the following two theorems.

\begin{Theorem}\label{sigmainSn}
Let $G = A \rtimes \mathbb S_n$ with $n>2$.  If {\rm dim}
$\mathfrak B(\mathcal{O}_{\sigma}^G, \rho) < \infty $, then  some of
the following hold with $Y_i \subseteq W_\chi$ or with $Y_i \cap W_\chi = \emptyset$:

\renewcommand{\theenumi}{\roman{enumi}}   \renewcommand{\labelenumi}{(\theenumi)}

\begin{enumerate}
\item $(1^{\lambda_1}, 2)$, $\mu_1 ={\rm sgn}$ or $\epsilon$, $\mu_2 =\chi _{(1; 2)}$.
\item $(2, 3)$, $\mu_2 =  \chi _{(1; 2)}$, $\mu_3 =\chi _{(0;3)}$.
\item $(  2^3)$, $\mu_2 =
\chi_{(1, 1, 1;2)}\otimes \epsilon$ or $\chi_{(1,1,1;2)}\otimes {\rm sgn}$.
\item  (4),  $\mu_4 =  \chi _{(2; 4)}$.
\end{enumerate}
\end{Theorem}

\begin{Theorem}\label{reducible}  Let $G= A\rtimes
\mathbb{S}_n$ be a classical Weyl group with $A\subseteq  (\mathbb Z_2)^n $ and $n>2$.

{\rm (i)} Assume that  there exist $i\not= j$  such that $\sigma
_i$, $\sigma _j $ $\notin A$. If  $ {\rm dim }\mathfrak B
(M)<\infty$, then  $n = 4$;  the type of $\sigma _p$ is $(2^2)$ and
the sign of $\sigma _p$ is stable; $\sigma $ has a negative cycle when  $\sigma _p \notin A.$

{\rm (ii)} Assume that  there exists  $\sigma _i$ $\in A$. If  $
{\rm dim }\mathfrak B (M)<\infty$, then  $\sigma _i = (a, \tau)$
with $\tau ^2 =1.$

{\rm (iii)} Assume    $\sigma _i $ $\in A$ for $1\le i \le m$. If
$ {\rm dim }\mathfrak B (M)<\infty$, then  there  is
 at most  one $\sigma_i \not= (1, 1, \cdots, 1)$.

 {\rm (iv)}  If $\sigma _i= \alpha =: (1,
1,\cdots, 1) \in G$ and $ \rho ^{(i)} = \theta _{\chi ^{(i)}, \mu ^{(i)}}
=:  (\chi ^{(i)} \otimes \mu ^{(i)}) \uparrow _{G_{\chi{(i)}}^\alpha  }^{G^\alpha }\in  \widehat
  {G ^{\alpha}}$, $1\le i \le m$. Then    $\mathfrak B (M)$ is  finite
dimensional if and only if $\chi ^{(i)}(\alpha) = -1$ for $i=1, 2,\cdots, m$.

\end {Theorem}

We prove Theorem \ref {sigmainSn} in Subsection \ref {s1.3}.
Theorem \ref {reducible}({\rm i}) follows from Remark \ref {3.12'}.
Theorem \ref {reducible}({\rm ii}) follows from Theorem  \ref {3.8}.
Theorem \ref {reducible}({\rm iii}) follows from Theorem  \ref {3.8}.
Theorem \ref {reducible}({\rm iv}) follows from Remark \ref {3.8'}.

Table 1 below lists the cases treated in this paper. In this table Nichols algebras  $ \mathfrak B({\mathcal O}_\sigma^G, \rho) $
of irreducible Yetter-Drinfeld modules over $G= A\rtimes \mathbb{S}_n$ ( $n>2$) are  infinite dimensional.

\begin {table}
\begin{tabular}{|l|l|l|l|}
  \hline
  { \bf case } & ${\bf \sigma}$ & {\bf Representation }& {\bf Reference } \\
\hline
  1&$  \sigma =e $, the unity of $G$& any  & See e.g. \cite {AZ07}   \\
\hline
  2&$  \sigma \in \mathbb A_n, n\ge 5,  $& any  & Theorem \ref {2.4}    \\
\hline
  3 &$  \sigma \in \mathbb S_n, n\ge 3, n\not=4, 5,6,$ & & \\
  &  $\sigma$ is not a transposition  & any  & Theorem \ref{sigmainSn}   \\
\hline
 4&the type of  $\sigma $ is $(2^3)$, $\sigma \in \mathbb S_6$, $ n=6$ &
 $\mu_2 \not=\chi_{(1, 1, 1;2)}\otimes \epsilon$, & \\
  & & $\mu_2 \not= $ $\chi_{(1,1,1;2)}\otimes {\rm sgn}$. & Theorem \ref{sigmainSn}\\
\hline
 5& the type of  $\sigma $ is $(2, 3)$, $\sigma \in \mathbb S_5$, $
n=5$ & $\mu_2 \not=  \chi _{(1; 2)}$, $\mu_3 \not=\chi _{(0;3)}$ & Theorem \ref{sigmainSn} \\
\hline
6& the type of  $\sigma $ is $4$, $\sigma \in \mathbb S_4$, $ n=4$
& $\mu_4 \not=  \chi _{(2; 4)}$ & Theorem \ref{sigmainSn} \\
\hline
 7&the type of  $\sigma $ is not  $ (2^3)$, $(2, 3)$,  $(4 )$, && \\
 &$(1^{\lambda_1}, 2)$ & any & Theorem \ref{sigmainSn}   \\
\hline
8&    $\sigma = (\alpha, \tau) $, $\alpha = (1, 1, \cdots, 1)$, &$\chi (\alpha) =1$,  &  \\
 & $\tau \in \mathbb S_n$, $\tau$ as in cases 1--7 &  $\mu$ as in cases 1--7 & Proposition \ref {2.6}\\
\hline
9&    $\sigma = (\alpha, \tau) $, $\alpha = (1, 1, \cdots, 1)$, & &  \\
 & $\tau \in \mathbb S_n$, the type of $\tau$ is not  &&\\
  &in list in Theorem \ref {2.6'} &  any  & Theorem \ref {2.6'}\\
\hline
10 &  $\sigma = (b, \tau)$  can be decomposed &&\\
& the multiplication of
 independent & &\\
& positive cycles.
the type of  $\tau  $ is not &&\\
 & $ (2^3)$, $(2, 3)$,  $(4 )$, $(1^{\lambda_1}, 2)$ & any & Theorem \ref {2.7'}   \\
\hline
11 &  $(b-\alpha, \tau)$  can be decomposed & &\\
& the multiplication of
 independent & &\\
& positive cycles.
the type of  $\tau  $ is not & &\\
 & in list in Theorem \ref {2.6'}. $\sigma = (b, \tau)$. & any & Theorem \ref {2.6'}   \\
\hline
 \end{tabular}
$$\hbox {Table } 1$$
\end {table}


This paper is organized as follows.
In section \ref{s2} we prove that except in three  cases Nichols
algebras of irreducible {\rm YD} modules supported by $\mathbb S_n $ are infinite dimensional.
In section \ref{s3}, we give a necessary and sufficient condition  for  a
Nichols algebra of a {\rm YD} module supported by  $A$ to be finite
dimensional. It is proved that  if $M$ is a reducible {\rm YD}
module over $ {\mathbb C } G$ supported by $\mathbb S_n$ with $n \ge 3,$ then
${\rm dim } {\mathfrak B } (M) = \infty$ and if $M$ is a  {\rm YD}
module over $ {\mathbb C } G$ supported by $\mathbb A_n$ with $n \ge 5$ ,  then  ${\rm dim } {\mathfrak B } (M) = \infty$. In
section \ref{s4}  we establish the relationship between Nichols algebras
over the Weyl groups of $B_n$ and $D_n$.
In the Appendix, the conjugacy classes of the Weyl groups of $B_n$ and $D_n$ are presented.

\section*{Preliminaries and Conventions}

Let ${\mathbb C}$ be the complex field and $G$ a finite group. Let  $\widehat{{
G}}$ denote the set of all isomorphism classes of irreducible
representations of the group $G$, $G^\sigma$ be the centralizer of
$\sigma$,    ${\mathcal O}_{\sigma}$ or ${\mathcal O}_{\sigma}^G$ be
the conjugacy class of $\sigma$ in $G$. If $D$ is a subgroup of $G$ and $C$ is a conjugacy
class of $D$, then $C_G$ denotes the conjugacy class of $G$ containing $C$.
The Weyl groups of $A_n$, $B_n$, $C_n$ and $D_n$ are
called the classical Weyl groups, written as $W(A_n)$, $W(B_n)$, $W(C_n)$ and $W(D_n)$, respectively.
Given a representation $\rho$ of the subgroup $D$ of $G$, let $ \rho \uparrow _D^G$ denote the induced
representation of $G$ as in \cite {Sa01}.

Let ${\mathcal K}(G)$ denote the set of
conjugacy classes $C$ in $G$.

For  $s\in G$ and  $(\rho, V) \in  \widehat {G^s}$, here is a
precise description of the {\rm YD} module $M({\mathcal O}_s,
\rho)$, see \cite {Gr00,AZ07}. Let $t_1 = s$, \dots, $t_{m}$ be a
numeration of  the conjugacy classes ${\mathcal O}_s$
containing $s$,  and choose $h_i\in G$ such that $h_i \rhd s =: h_i s
h_i^{-1} = t_i$ for all $1\le i \le m$. Then $M({\mathcal O}_s,
\rho) = \oplus_{1\le i \le m}h_i\otimes V$. Let $h_iv =: h_i\otimes
v \in M({\mathcal O}_s,\rho)$, $1\le i \le m$, $v\in V$. If $v\in V$
and $1\le i \le m $, then the action of $h\in G$ and the coaction are given by
\begin {eqnarray} \label {e0.11}
\delta(h_iv) = t_i\otimes h_iv, \qquad h\cdot (h_iv) =h_j(\gamma\cdot v),
\end {eqnarray}
 where $hh_i = h_j\gamma$, for some $1\le j \le m$ and $\gamma\in G^s$.  Let $\mathfrak{B}
({\mathcal O}_s, \rho )$ denote $\mathfrak{B} (M ({\mathcal O}_s,
\rho ))$.    A {\rm YD} module is said to be reducible if it
has a non-trivial {\rm YD} submodule. If $\rho = \rho ^{(1)} \oplus
\rho ^{(2)} \oplus \cdots \oplus \rho ^{(m)}$ with $\rho ^{(i)} \in
\widehat {G^{\sigma}}$ for $1\le i\le m$, then $ M(\mathcal
{O}_\sigma ^G, \rho ^{(1)}) \oplus M(\mathcal {O}_\sigma ^G, \rho
^{(2)}) \oplus \cdots \oplus M(\mathcal {O}_\sigma ^G, \rho ^{(m)})$
is called a {\rm YD} module of $\rho$, also written as $M(\mathcal
O_\sigma^G, \rho)$.

Let $V$ be a braided vector space of diagonal type with a basis $\{x_i \mid i= 1, 2,
\cdots, n\}$ and $B(x_i\otimes x_j) = q_{ij}(x_j\otimes x_i)$. If
there exists a generalized Cartan matrix $(a_{ij})$ such that     satisfy
\begin {eqnarray} \label {e0.0.1}
q_{ij}q_{ji} = q_{jj} ^{a_{ji}}
\end {eqnarray}
for any $i, j= 1, 2,\cdots, n$, then braiding $B$ (or $V$, or $\mathfrak B(V)$) is
called a braiding of the Cartan type.
We assume that we choose $a_{ij}$ such that they satisfy (\ref {e0.0.1}).
That is, $a_{ij}$ is the maximal non-positive integer satisfying (\ref {e0.0.1}).  Thus
$\mathfrak B(V)$ is finite dimensional if and only if $A$ is of
finite type (see Theorem 4 of \cite  {He06}).

\section {Classical Weyl groups} \label {s2}

In this section we give a necessary and sufficient
condition  for a Nichols algebra of irreducible {\rm YD} module
supported by $A$ to be finite dimensional, and show that  except in three  cases Nichols algebras
of irreducible {\rm YD} modules supported by $\mathbb S_n $ are infinite dimensional.

Let $H= B \rtimes  D$ be a semidirect product of $B$ and $D$, where $B$ is abelian groups.
 For any $\chi \in \widehat  B,$ let $D_\chi =: \{ h \in D \mid h\cdot \chi =
\chi\}$; $H_\chi =: B \rtimes D_\chi$. For an irreducible
representation $\rho $ of $D_\chi$, let  \begin {eqnarray} \label {e0.0.5}\theta _{\chi, \rho} =:
(\chi \otimes \rho) \uparrow ^H_{H_\chi}\end {eqnarray} denote the induced
representation of $\chi \otimes \rho$ on $H$. By Proposition 2.5 of \cite
{Se77}, every irreducible representation of $H$ is of the
 form: $\theta _{\chi, \rho}$. Let $\epsilon \in \widehat  B$
with $\epsilon (a) =1$ for any $a \in B.$ Thus $D_\epsilon = D$ and
$\theta _{\epsilon, \rho}$ is an irreducible representation of $H$.

\begin {Lemma}\label {1.1'} Let  $H = B \rtimes D$.

 {\rm (i)} If  $\sigma \in D,$ then $H^\sigma =B^\sigma \rtimes D^\sigma.$

 {\rm (ii)} $\theta _{\chi, \rho}$ ( as in (\ref {e0.0.5}) ) is a one
dimensional representation of $H^\sigma= B^\sigma \rtimes D^\sigma$
if and only if $D_\chi^\sigma =D^\sigma$ and ${\rm deg} \rho =1$.

{\rm (iii)} If $\sigma \in B$,
then $(B \rtimes D)^\sigma = B \rtimes D^\sigma$.

\end {Lemma}
\noindent {\bf Proof.} {\rm (i)} If $x = (a, d) \in H^\sigma,$ then $x\sigma = \sigma
x.$ Thus \begin {eqnarray} \label {e1.11} a  = \sigma \cdot a \ \
\hbox {and } &&   d\sigma = \sigma d. \end {eqnarray} This implies
$d\in D^\sigma$ and $a\in B^\sigma $ since $\sigma \cdot a = \sigma
a \sigma ^{-1}$.

Conversely, if $x = (a, d) \in B^\sigma \rtimes D^\sigma$, then
(\ref {e1.11}) holds. This implies $x \sigma = \sigma x$ and $x \in
G^\sigma.$

 {\rm (ii)} Let $P$ and $X$ be the representation spaces of $\chi$
and $\rho$ on $B^\sigma$ and $D^\sigma_\chi$, respectively. $ (( P
\otimes X) \otimes_{{\mathbb C } G_\chi^\sigma} {\mathbb C } H^\sigma , \theta _{\chi,
\rho})$ is a one dimensional representation of $H^\sigma= B^\sigma
\rtimes D^\sigma$ if and only if  $ {\mathbb C }  H^\sigma ={\mathbb C }  H^\sigma_\chi$ and
${\rm dim } X= 1$. However. $ {\mathbb C }  H^\sigma =  {\mathbb C }  H^\sigma_\chi$ if and
only if $D_\chi ^\sigma = D^\sigma$.

{\rm (iii)} It is clear that $(a, \tau) \in (B \rtimes D)^\sigma $ if and only if $\tau \in D^\sigma$.
$\Box$

\vskip.1in
Consequently, $\theta _{\chi, \rho}= \chi  \otimes \rho$ when
$\theta _{\chi, \rho}$ is one dimensional representation.

\subsection {$\sigma \in A$} \label {s1.1}

\vskip.1in
Let $\sigma = ({a_1}, {a_2}, \cdots, {a_2} )\in A$.  If
$\rho $ is  a  representation of $\widehat{G^\sigma}$, then $\rho =\theta _{
(\chi \otimes \mu)} = (\chi\otimes \mu ) \uparrow
^{G^\sigma}_{(G^\sigma)_\chi}$ with $ \chi = \chi_2^{b_1} \otimes
\chi_2^{b_2} \otimes \cdots \otimes \chi_2^{b_n} $ $\in \widehat  A$,
$(G^\sigma)_\chi = A \rtimes (\mathbb S_n ^\sigma)_\chi$ and
$(\mathbb S_n ^\sigma)_\chi =: $ $\{ \tau \in \mathbb S_n^\sigma
\mid \tau \cdot \chi = \chi \}$ and
$\mu \in \widehat { (\mathbb S_n^\sigma )_{\chi}}$, see Proposition 2.5 of \cite
{Se77}. Let $f_{\tau\cdot \sigma, \chi } =:
{b_1}{a_{\tau^{-1}(1)}} + {b_2}{a_{\tau^{-1}(2)}} + \ldots
{b_n}{a_{\tau^{-1}(n)}}$, $W_\sigma =: \{i \mid  a_i \neq 0 \}$ and
$W_\chi  = \{i \mid  b_i \neq 0 \}$. It is clear that $f _{\tau
\cdot \sigma, \chi } = \mid W_{\tau \cdot \sigma} \cap W_\chi \mid
$, where $W_{\tau \cdot \sigma} = \{ i \mid a_{\tau ^{-1} (i)} \not=
0 \} = \tau W_{\sigma}.$ $f_{\tau\cdot \sigma, \chi }$ is written as
$f_{\tau\cdot \sigma}$ in short. Note  that $\mathbb{S}_n$ acts on
$\widehat  A$ as follows: for any $\chi \in \widehat  A$ with $\chi  = \chi
_2^{b_1} \otimes \chi _2^{b_2} \otimes \cdots  \chi _2^{b_n} $ and
$h\in \mathbb{S}_n$, define
 \begin {eqnarray*} h \cdot \chi =:
( \chi_2^{b_{h^{-1}(1)}}, \chi_2^{b_{h^{-1}(2)}}, \cdots,
\chi_2^{b_{h^{-1}(n)}} ).\end {eqnarray*}

\begin {Lemma}\label {2.2} Under the notations above, we have

{\rm (i)} $\rho (\sigma)= q_{\sigma, \sigma} {\rm id}$ with
$q_{\sigma, \sigma} = (-1) ^{f_\sigma}$.

{\rm (ii)} If  $W_\sigma = W_\chi$, then $(\mathbb S_n^\sigma)_\chi
=\mathbb S_n^\sigma $.

{\rm (iii)} If  $W_\sigma = W_\chi$, then $f_{g\cdot \sigma} =
f_{g^{-1}\cdot  \sigma }$ for any $g\in \mathbb S_n$.

\end {Lemma}

\noindent {\bf Proof.}  Let $P$ and $V$ be the representation spaces of $\chi$
and $\mu$, respectively. Then the  representation space of $\rho$ is
$ {\mathbb C } G^\sigma \otimes _{k(G^\sigma )_\chi } (P \otimes V)$.  For any $p\in P$, $v\in V$ and
$\tau \in \mathbb S_n$, we have
\begin {eqnarray}\label {e2.2.1} &&\rho (\tau \cdot \sigma) (1 \otimes_{k(G^\sigma )_\chi }
(p \otimes v)) =  1 \otimes_{k(G^\sigma )_\chi } (( (\tau\cdot \sigma)\cdot p) \otimes x )\\
&=& \chi (\tau \cdot \sigma) 1 \otimes_{k(G^\sigma )_\chi } (p
\otimes x) = (-1) ^ {f_{\tau \cdot \sigma}}  1 \otimes_{k(G^\sigma
)_\chi } (p \otimes x). \nonumber\end {eqnarray}

{\rm (i)}  It follows from (\ref {e2.2.1}) that  $q_{\sigma, \sigma
} = \chi (\sigma)= (-1)^{f_\sigma}. $

{\rm (ii)} If $\tau  \in \mathbb S_n^\sigma,$ then $\tau ( W_\sigma
) = W_\sigma$, i.e. $\tau( W_\chi ) = W_\chi$. Consequently, $\tau
\in (\mathbb S_n^\sigma)_\chi $.

{\rm (iii)} Let $g(W_\sigma) \cap W_\sigma= \{l_1, l_2, \cdots,
l_r\}$ with $g(m_i) = l_i$ for $1\le i\le r.$ Then $g^{-1}(l_i) =
m_i$ for $1\le i\le r$ and  $g ^{-1}(W_\sigma) \cap W_\sigma  =
\{m_1, m_2, \cdots, m_r\}$. This implies $f_{g\cdot \sigma} =
f_{g^{-1}\cdot \sigma }$.
 $\Box$

\begin {Theorem}\label {2.3} Let $G = A \rtimes \mathbb S_n$ with $\sigma \in A$ and
$n>2$. Then {\rm dim}
$\mathfrak B(\mathcal{O}_{\sigma}^G,\rho) < \infty $ if and only if
$f_\sigma$ is odd and either $W_\chi = W_\sigma$ or $|W_\sigma|=n$ or $|W_\chi|=n$.
 \end {Theorem}

 \noindent {\bf Proof.}  Let $P$ and $V$ be the representation spaces of $\chi$
and $\mu$, respectively. Thus  the representation space of $\rho$ is
$ {\mathbb C } G^\sigma \otimes _{k(G^\sigma )_\chi } (P \otimes V)$. Let $\xi =
1 \otimes _{k(G^\sigma )_\chi } (p_0 \otimes v_0) $ with $0\not=
p_0\in P$ and $0\not= v_0 \in V$.

We show this by the following seven steps.

 {\rm (i)} If $f_\sigma$ is even, then $\rho (\sigma) = (-1) ^{f_{\sigma}} {\rm id} = {\rm id}$
 by Lemma \ref {2.2} (see also Remark 1.1 in \cite{AZ07}) and {\rm dim}
$\mathfrak B(\mathcal{O}_{\sigma}^G,\rho) = \infty $.

From now on we assume that $f_\sigma $ is odd.

{\rm (ii) } Assume that $i_0 \in W_\chi \cap W_\sigma$, $i_1 \in
W_\sigma \setminus W_\chi$ and $i_2 \notin W_\sigma$. Let $h_1 =1$,
$h_2 = (i_0, i_1, i_2)$, $h_3 = h_2 ^{-1}$, $t_i = h_i \cdot \sigma$
and $\gamma _{ij} = h_j ^{-1}h_i \cdot \sigma$ for $i, j = 1, 2, 3.$
By simple computation we have $\gamma_{12} = h_2^{-1}\cdot \sigma, $
$\gamma_{21} = h_2\cdot \sigma, $ $\gamma_{13} = h_2\cdot \sigma, $
$\gamma_{31} = h_2^{-1}\cdot \sigma, $ $\gamma_{23} = h_2^{-1}\cdot
\sigma, $ and $\gamma_{32} = h_2\cdot \sigma.$
 Since $(h_2\cdot \sigma ) _{i_0} = 0$, $i_0$-th component
 of $h_2\cdot \sigma$. $(h_2\cdot \sigma ) _{i_1} = 1$,  $(h_2^{-1}\cdot \sigma )
_{i_0} = 1$, and $(h_2^{-1}\cdot \sigma ) _{i_1} = 0$, we
have that $\sigma$, $h_2\cdot \sigma$ and $h_2^{-1}\cdot \sigma$ are not the same. Indeed,
\begin {eqnarray*}
f_{h_2\cdot \sigma } +f_{h_2^{-1}\cdot \sigma }
&\equiv & (a_{h_2^{-1}(i_0)} + a_{h_2(i_0)}) b_{i_0} +
(a_{h_2^{-1}(i_1)} + a_{h_2(i_1)}) b_{i_1} + (a_{h_2^{-1}(i_2)} +
a_{h_2(i_2)}) b_{i_2} \\
&\equiv& (a_{i_2} + a_{i_1}) b_{i_0} + (a_{i_1} + a_{i_0})b_{i_2}
\equiv 1 \ \ (\hbox {mod} \   2).
\end {eqnarray*}
Consequently, $f_{h_j ^{-1}h_i\cdot \sigma } +f_{h_i ^{-1}h_j \cdot \sigma }
\equiv 1 \ \ (  \hbox { mod } \  2)$ for $i\not= j.$ Let $B(h_i \xi
\otimes h_j\xi) = q_{ij} (h_j \xi \otimes h_i\xi)$, where $B$
is the braiding of the Cartan type. Then $q_{ij} = (-1)^{f_{h_j^{-1}h_i\cdot
\sigma}}$ by (\ref {e2.2.1}) and $q_{ij} q_{ji}=-1 = q_{ii}
^{a_{ij}}$ with $a_{ij} = -1$ when $i\not= j.$ The braided subspace
spanned by $\{ h_i \xi \mid i =1, 2, 3\}$ is not of finite Cartan type.
Consequently, {\rm dim} $\mathfrak B(\mathcal{O}_{\sigma}^G,\rho) =
\infty $ by Theorem 4 of \cite  {He06}.

{\rm (iii)} Assume that  $i_0 \in W_\chi \cap W_\sigma$, $i_1 \in
W_\chi\setminus W_\sigma $ and $i_2 \notin W_\chi$. Let $h_1 =1$,
$h_2 = (i_0, i_1, i_2)$, $h_3 = h_2 ^{-1}$,  $t_i = h_i \cdot
\sigma$ and $\gamma _{ij} = h_j ^{-1}h_i \cdot \sigma$ for $i, j =
1, 2, 3.$
 Since $(h_2\cdot \sigma ) _{i_1} = 1$,
$(h_2^{-1}\cdot \sigma ) _{i_0} = 0$,  $(h_2^{-1}\cdot \sigma )
_{i_2} = 1$, and $(h_2\cdot \sigma ) _{i_2} = g_2^0$, we have
that $\sigma$, $h_2\cdot \sigma$ and $h_2^{-1}\cdot \sigma$ are not
the same. See
\begin {eqnarray*}
f_{h_2\cdot \sigma } +f_{h_2^{-1}\cdot \sigma }
&\equiv & (a_{i_2} + a_{i_1}) b_{i_0} + (a_{i_0}+ a_{i_2}) b_{i_1}
\equiv 1 \ \ ( \hbox {mod} \  2).
\end {eqnarray*}
Consequently, {\rm dim} $\mathfrak
B(\mathcal{O}_{\sigma}^G,\rho) = \infty $ by  the same arguments as in Part (ii).

{\rm (iv)} If $\mid W_\sigma \mid =n$, then {\rm dim} $\mathfrak
B(\mathcal{O}_{\sigma}^G,\rho) < \infty $ since it is a central quantum linear space.

{\rm (v)}  Assume $W_\chi = W_\sigma$. Then $(\mathbb S_n ^\sigma )_\chi = \mathbb S_n ^\sigma
 $  by Lemma \ref {2.2}. Let $v_1, v_2, \cdots v_s$ be the basis of $V$ and $\xi _i =
  1 \otimes _{k(G^\sigma )_\chi } (p_0 \otimes v_i) $ for $1\le i \le
  s$. $\{h_i \xi _j \mid 1\le i \le m, 1\le j \le s \}$ is a basis
  of  $M (\mathcal O_\sigma ^G, \rho)$, where $G = \cup _{i=1}^m h_i G^\sigma
  $ is a left coset decomposition of $G^\sigma$ with $h_i \in \mathbb S_n.$
  Let $t_i =: h_i \cdot \sigma = h_i \sigma h_i^{-1}$ and $\gamma _{ij}  =: h_j ^{-1} h_i \cdot
   \sigma $. Therefore, $t_ih_j = h_j r_{ij}$ and $B(h_i \xi _u \otimes h_j \xi _v )=
   \chi (h_j^{-1}h_i\cdot \sigma  ) (h_j\xi _v \otimes h_i \xi _u)$. By Lemma \ref {2.2}
   $\chi (h_j^{-1}h_i\cdot \sigma ) \chi (h_i^{-1}h_j \cdot \sigma  ) =1$
    for any $1\le i, j \le    m$, $1\le u, v \le
  s.$ Consequently, $M({\mathcal O } _\sigma ^G, \rho)$ is of a finite
  Cartan type and {\rm dim} $\mathfrak B(\mathcal{O}_{\sigma}^G,\rho) < \infty.$

{\rm (vi)} Assume $\mid  W_\chi \cap W_\sigma \mid \ge 1$ with $\mid
W_\chi \mid =n$. It is clear $(\mathbb S_n^\sigma  )_\chi = \mathbb
S_n ^\sigma.$ Since $f_{\sigma}$ is odd, $f_{\tau\cdot \sigma}$ is
odd for any $\tau \in \mathbb S_n$. As the proof of (v), {\rm dim}
$\mathfrak B(\mathcal{O}_{\sigma}^G,\rho) < \infty.$

{\rm  (vii)}  If {\rm dim} $\mathfrak B(\mathcal{O}_{\sigma}^G,\rho)
< \infty,$ then $f_\sigma$ is odd. If $W_\sigma \not= W_\chi$, then
$\mid W_\sigma \setminus W_\chi \mid \ge 1$ or $\mid W_\chi
\setminus W_\sigma \mid \ge 1$. When  $\mid W_\sigma \setminus
W_\chi \mid \ge 1$ with $\mid W_\sigma \mid < n$, it is  case {\rm
(ii)}. When $\mid W_\chi \setminus W_\sigma \mid \ge 1$ and $\mid
W_\chi \mid <n$ it is case {\rm (iii)}. When $\mid W_\sigma \mid =n$
it is case {\rm (iv)}. When $\mid  W_\chi \cap W_\sigma \mid \ge 1$
with $\mid W_\chi \mid =n$, it is case {\rm (vi)}.
 $\Box$

\begin {Corollary}\label {2.3'} Let $G = A \rtimes \mathbb S_n$ with $\sigma \in A$ and
$n>2$. Assume $\rho = \theta _{\chi, \mu} = (\chi\otimes \mu )
\uparrow ^{G^\sigma}_{(G^\sigma)_\chi}$. If
 {\rm dim} $\mathfrak B(\mathcal{O}_{\sigma}^G,\rho) < \infty $, then
\begin {eqnarray} \label {e2.3.1}
\mid g(W_{\sigma }) \cap W_{\chi} \mid \equiv
  \mid g^{-1}(W_{\sigma }) \cap W_{\chi} \mid \ \ ( \hbox {  mod  } \ 2)
  \end {eqnarray}
  for any $g \in \mathbb S_n $.
 \end {Corollary}

 \noindent {\bf Proof.}  Suppose that there exists $g \in \mathbb S_n $
 such that (\ref {e2.3.1}) does not hold. This implies that $f_{g\cdot \sigma, \chi} \equiv
f_{g\cdot \sigma, \chi}  \ (  {\rm mod } \ 2) $ does not hold. By
Lemma \ref {2.2}, $W_\sigma \not= W_{\chi}. $ Considering {\rm dim}
$\mathfrak B(\mathcal{O}_{\sigma}^G,\rho) < \infty $ and Theorem
\ref {2.3}, we have that $\mid W_\chi \mid = n$. See
\begin {eqnarray*} \mid g(W_{\sigma }) \cap W_{\chi} \mid = \mid
g W_\sigma  \mid = \mid g^{-1} W_\sigma  \mid =
  \mid g^{-1}(W_{\sigma }) \cap W_{\chi} \mid. \
\end {eqnarray*} This is a contradiction. $\Box$

\subsection {$\sigma \in \mathbb A_n, $ alternating group}\label {s1.2}

\begin {Theorem}\label {2.4} Let $G = A \rtimes \mathbb{S}_n$.
If $\sigma \in \mathbb A_n $ with $n\ge 5$, then
 {\rm dim} $\mathfrak
B(\mathcal{O}_{\sigma}^G, \rho) = \infty $ for any $\rho \in
\widehat {G^\sigma}$.
\end {Theorem}
\noindent {\bf Proof.} It follows from Theorem 1.2 of \cite {AFGV} that {\rm
dim} $\mathfrak B(\mathcal{O}_{\sigma}^{\mathbb A_n}, \rho \mid
_{\mathbb A_n ^\sigma}) = \infty $. Applying Lemma 3.2 in \cite 
{AFGV}, we have {\rm dim} $\mathfrak B(\mathcal{O}_{\sigma}^G,
\rho) = \infty$. $\Box$

\subsection {$\sigma \in \mathbb S_n$} \label {s1.3}

We rely on the general theory of representations of $\mathbb S_n$ as described on Pages 295-299 of \cite{Su78}.
Let $\sigma\in \mathbb S_n$  with cycle  type $(1^{\lambda _1}, 2^{\lambda _2}, \cdots,
n^{\lambda _n})$.  Let  $r_j=:\sum_{1\leq k
\leq j-1}k\lambda_k$ and $\sigma _j =: \prod _{1\le l \le \lambda
_j}$ $\Big( y_{ r_j+(l-1)j+1}, $  $ \qquad y_{r_j+(l-1)j+2}, \quad
\cdots, \quad y_{ r_j+lj }\Big)$, the multiplication of cycles of
length $j$ in the independent cycle decomposition of $\sigma$, as
well as  $Y_j =: \{y_s \mid s=
r_{j}+1, \cdots, r_{j+1}\}$.  Therefore $\sigma= \prod \sigma_i$ and
$ ({\mathbb S}_n)^\sigma= \prod  ({\mathbb S}_{Y_i})^{\sigma_i}$
$=T_1 \times \cdots \times T_n.$ Obviously, every element in
$\widehat {(C_l)^m}$ can be denoted by $\chi _{(t_{1,l}, t_{2, l},
\cdots, t_{m,l}; l)}=: \chi _l ^{t_{1,l}} \otimes \chi _l
^{t_{2,l}}\otimes \cdots \otimes \chi _l ^{t_{m,l}}$ for  $0 \le
t_{j, l}\le l-1$.

\begin {Proposition}\label {3.7'} Let $G= A \rtimes  \mathbb{S}_n$ with
  $\sigma \in \mathbb{S}_n$. Then
  $(\mathbb S_n)^\sigma_{\chi^{(\nu)}} = (\mathbb S_n)^\sigma $ if and only if
  $Y _j \subseteq \{ 1, 2, \cdots, \nu \}$ or  $Y _j \subseteq  \{\nu +1, \nu+2, \cdots, n\}$
   for $1\le j\le n $, where $Y_j$ is the same as in  the beginning of Section \ref {s2}.
\end {Proposition}
\noindent {\bf Proof.}
  If $Y _j \subseteq \{ 1, 2, \cdots, \nu \}$, then $A_{l, j},
B_{h, j} \in {\mathbb S}_{\{1, 2, \cdots, \nu\}} $ for $1\le l \le
\lambda _j$ and $1\le h\le \lambda _j-1$, where $A_{l, j}$, $B_{h,
j}$ are the same as in  the begin of Section \ref {s2}. If $Y _j
\subseteq \{ \nu +1, \nu +2, \cdots, n \}$, then $A_{l, j}, B_{h, j}
\in {\mathbb S}_{\{\nu +1, \nu +2, \cdots, n\}} $ for $1\le l \le
\lambda _j$ and $1\le h\le \lambda _j-1$. Consequently, $(\mathbb
S_n)^\sigma_{\chi^{(\nu)}} = (\mathbb S_n)^\sigma$

Conversely, assume $(\mathbb S_n)^\sigma_{\chi^{(\nu)}} = (\mathbb
S_n)^\sigma $. If there exists $1\le j \le n$
 such that  $Y _j \nsubseteq \{ 1, 2, \cdots, \nu \}$ and $Y _j \nsubseteq  \{\nu +1, \nu+2, \cdots,
 n\}$, then there exist $a, b \in Y_j$ with $a\in \{ 1, 2, \cdots, \nu
 \}$ and $b \in \{\nu +1, \nu+2, \cdots,
 n\} $. Note  $Y_j = \{y_s \mid s=r_{j}+1, \cdots, r_{j+1}\}$.

If there exists $l$ such that $a, b\in \{  y_{ r_j+(l-1)j+1}, \qquad
y_{r_j+(l-1)j+2},  \quad \cdots,  \quad y_{ r_j+lj }\}$, then $A_{l,
j} \notin (\mathbb S_n)^\sigma_{\chi^{(\nu)}}$, a contradiction.
Thus there exist $l\not=l'$ such that $a\in A_{l, j}$ and $b \in
A_{l', j}$. Let $a = y_{ r_j+(l-1)j+s} $ and $b=  y_{
r_j+(l'-1)j+s'} $. Considering $B_{h,j} \in (\mathbb
S_n)^\sigma_{\chi^{(\nu)}}$, we have $y_{ r_j+(l'-1)j+s} \in \{ 1,
2, \cdots, \nu  \}$ which is a contradiction.
$\Box$

 Let $\sigma\in \mathbb
S_n$ be of type $(1^{\lambda_1}, 2^{\lambda_2}, \dots,
n^{\lambda_n})$. Assume that   $\rho = \theta _{\chi , \mu } = (\chi
\otimes \mu )\uparrow ^{G^\sigma} _{(G^\sigma )_\chi} \in \widehat
{A^\sigma \rtimes \mathbb S_n ^\sigma }$ with $\chi = (\chi _2
^{b_1} \otimes \cdots \chi_2 ^{b_n})  \in \widehat {A^{\sigma}}$ and
$\mu \in \widehat {({\mathbb S}_n^\sigma )_\chi}$.


\renewcommand{\theenumi}{\roman{enumi}}   \renewcommand{\labelenumi}{(\theenumi)}


\vskip.1in
We now in the position to prove the statements in Theorem \ref{sigmainSn}. 

{\bf Proof of Theorem \ref{sigmainSn}.} Let $\rho \mid _{\mathbb S_n^\sigma } = \rho ^{(1)}
\oplus\rho ^{(2)} \oplus\cdots \oplus \rho ^{(m)}$ be a
decomposition of simple modules of $\rho\mid _{\mathbb S_n^\sigma
}$. Since {\rm dim} $\mathfrak B(\mathcal{O}_{\sigma}^G, \rho) <
\infty $, we have  {\rm dim} $\mathfrak
B(M(\mathcal{O}_{\sigma}^{\mathbb S_n}, \rho ^{(1)}) \oplus
M(\mathcal{O}_{\sigma}^{\mathbb S_n}, \rho ^{(2)}) \oplus \cdots
\oplus  M(\mathcal{O}_{\sigma}^{\mathbb S_n}, \rho ^{(m)}) )< \infty
$ by Lemma 2.2 of \cite {AFGV}. Considering Corollary
8.4 in \cite {HS10} we have $m=1$ and $\rho \mid _{\mathbb S_n^\sigma}$ is
irreducible. Furthermore, $\rho \mid _{\mathbb S_n^\sigma}$ is one
dimensional by  Theorem 1 of \cite {AFZ08} or Theorem 1.1 of \cite {AFGV}. 
Since ${\rm deg } (\rho) ={\rm deg } ( \rho \mid
_{\mathbb S_n^\sigma}) $, we have that $\rho$ is one dimensional and
$(\mathbb S_n^\sigma)_\chi = \mathbb S_n^\sigma$ by  Lemma \ref
{1.1'}. Consequently, $Y _j \subseteq W_\chi$ or $Y _j \cap W_\chi
=\emptyset$
 by Proposition \ref {3.7'}. Consequently, the result  follows from  Theorem 1.1 in \cite
 {AFGV} and Tabel 2 in \cite {AZ07}.
 $\Box$

\subsection { $\sigma = (b, \tau) $ with $ \sum _{i=1}^n b_i=0$}\label {s1.4}

\begin {Lemma}\label {2.7} (i) Assume   $b \in A$ and $\xi =: (1, 2, \cdots, n )\in {\mathbb S}_n$,  
then $\sum _{i=1}^n b_i \equiv 0$ if and only if there exists $a\in A$ such that $b \equiv (1-\xi) (a)$.

(ii) If $(b, \sigma) \in A\rtimes \mathbb S_n$, then $(b, \sigma)$ can be decomposed into a product of  
independent  positive cycles if and only if  $(b, \sigma )$ and $\sigma$ are conjugate to each other.

(iii) If $\alpha = (1, 1, \cdots, 1)\in A$  and $(b-\alpha , \sigma) \in A\rtimes \mathbb S_n$, 
then $(b-\alpha , \sigma)$ can be decomposed  into a product of  independent  positive cycles 
if and only if $(b, \sigma )$ and $(\alpha, \sigma)$ are conjugate to each other.

\end  {Lemma}
\noindent {\bf Proof.} (i) The sufficiency.  Set $ a_i = b_1 + \cdots + b_i$, for
$1\le i \le n.$  Thus $b \equiv (1-\xi) (a)$. The necessity is clear. 

(ii) The necessity.  Let $(b, \sigma ) = (b^{(1)}, \sigma _1)(b^{(2)}, \sigma _2)\cdots
(b^{(r)}, \sigma _r)$ is a decomposition of independent positive cycles. By part (i), there exists $a ^{(i)}$ such that
$b^{(i)} \equiv (1-\sigma_i) (a^{(i)})$ for $i= 1, 2, \cdots, r$. Let
$ a =: a ^{(1)} + \cdots + a^{(r)}$. It is clear $b \equiv (1-\sigma) (a)$ and
$ (a, \sigma) (0, \sigma) (a, \sigma) ^{-1} \equiv (b, \sigma)$.
The sufficiency follows from Theorem \ref {6.2} in the Appendix.

(iii) The necessity. By (ii), there exists $a\in A$ such that $(1-\sigma) (a) \equiv b-\alpha$ and 
$ (a, \sigma) (0, \sigma) (a, \sigma) ^{-1}\equiv (b-\alpha, \sigma)$. Considering $\alpha \in Z(G)$, 
we have $ (a, \sigma) (\alpha, \sigma) (a, \sigma) ^{-1} \equiv (b, \sigma)$. 
The sufficiency follows from Theorem \ref {6.2} in the Appendix.
$\Box$

\begin {Theorem}\label {2.7'} Assume that $(b, \sigma) \in A\rtimes \mathbb S_n$ can be decomposed into
a product of  independent  positive cycles, $n>2$. If the type of $\sigma$  is not
 $( 1^{\lambda _1}, 2)$ or  $(2, 3)$ or   $(2^3)$ or $(4)$, then {\rm dim}
$\mathfrak B(\mathcal{O}_{(b, \sigma)}^G, \rho) = \infty $ for any $\rho \in \widehat {G^{(b, \sigma)}}$.
\end {Theorem}

\noindent {\bf Proof.} By Theorem 1.1 in \cite {AFGV} and Theorem \ref{sigmainSn}, $\dim\mathfrak B(\mathcal{O}_{ \sigma}^G, \rho) = \infty $.  
Consequently, Our claim  follows from Lemma \ref {2.7}(ii). $\Box$

\subsection { $\sigma = (\alpha, \tau)$ with $\alpha = (1, 1, \cdots,1 ) \in A$}\label {s1.5}

Let $G = A \rtimes \mathbb S_n$ and $\tau \in \mathbb S_n$ with  $\sigma = (\alpha, \tau)
\in G.$ It is clear that the following holds:

{\rm (i)} $G^{\sigma} = G^{\tau} = A^\tau \rtimes \mathbb S_n ^\tau$.

{\rm (ii)} For any $\rho \in \widehat { G^\sigma }$, there exists
$\chi \in \widehat {A^\tau}$ and $\mu \in \widehat {(\mathbb S_n
^\tau)_\chi}$ such that $\rho = \theta _{\chi, \mu}= (\chi \otimes
\mu )\uparrow ^{G^\sigma} _{(G^\sigma )_\chi}$.

\begin {Proposition}\label {2.6} Under the notations above we have:

{\rm (i)} $\mathcal{O}_{(\alpha, \tau)}^G = \alpha \mathcal{O}_{(0, \tau)}^G$.

{\rm (ii)} If $H$ is a group and $\sigma \in H$, $a\in Z(H)$ (the center of $H$), 
then $\mathcal{O}_\sigma ^H$ is a rack  of type $D$ if and only if $a \mathcal{O}_\sigma ^H$ is a rack of type $D$ 
(which  was defined in \cite  {AFGV}, Definition 3.5).

{\rm (iii)} If $\chi (\alpha) =1$, then
 {\rm dim} $\mathfrak B(\mathcal{O}_{(\alpha, \tau)}^G, \rho) < \infty $ if
 and only if {\rm dim }$\mathfrak B(\mathcal{O}_{(0,  \tau)}^G, \rho ) < \infty $

{\rm (iv)} If $\chi (\alpha) =-1$, then it is impossible that
 $\mathfrak B(\mathcal{O}_{(\alpha, \tau)}^G, \rho)  $
and  $\mathfrak B(\mathcal{O}_{(0,  \tau)}^G, \rho )$ are  finite dimensional simultaneously.
\end {Proposition}

\noindent {\bf Proof.} {\rm (i)} It is clear.

(ii) It is clear that $X$ is a subrack of $\mathcal{O}_\sigma ^H$ if and only if
 $aX$ is a subrack of  $a \mathcal{O}_\sigma ^H$. If $\mathcal{O}_\sigma ^H$ is of type $D$, 
then there exists a decomposable subrack   $Y = R \cup S$ of $\mathcal{O}_\sigma ^H$. 
For any $r\in R, s\in S, $ $r \rhd (s \rhd (r \rhd s)) \not=s$ , that implies
 $ar \rhd (as \rhd (ar \rhd as)) \not=as$. Consequently, $aY = aR \cup aS$ of $a\mathcal{O}_\sigma ^H$ and 
$a \mathcal{O}_\sigma ^H$ is a rack of type $D$.  Conversely, the proof is similar.

{\rm (iii)}  Let $M=: M(\mathcal{O}_{(\alpha, \tau)}^G, \rho)$ and
$ N=: M(\mathcal{O}_{(0, \tau)}^G, \rho )$ with $\rho (\alpha, \tau)
= q _M \,{\rm id}$ and $\rho (0, \tau) = q _N {\rm id}.$ It is
clear that $q_M = \chi (\alpha) q_{N}$. Let $V$ be the
representation space of $\rho$ with a basis  $x_1, \cdots, x_n$.
$\mathcal{O}_{(\alpha, \tau)}^G = \alpha \mathcal{O}_{(0, \tau)}^G$.
Assume that $G= \cup _{i=1}^m h_iG^{(\alpha, \tau)}$ is the
decomposition of left cosets of $G^{(\alpha, \tau)}$ in $G$. Then
$G= \cup _{i=1}^m h_iG^{(0, \tau)}$ is the decomposition of left
cosets of $G^{(0, \tau)}$ in $G$. It is clear $\mathcal{O}_{(0,
\tau)}^G = \{t_i \mid  t_i = h_i \tau h_i ^{-1 }, 1\le i \le m \}$.
See
\begin {eqnarray}\label {e2.9}
B_M(h_i x_j \otimes h_s x_t) &=& \chi (\alpha ) B_N(h_i x_j \otimes
h_s x_t)
\end {eqnarray} for $1\le i, s\le m, $ $1\le j, t\le n,$ where $B_N$
and $B_M$ denote the braidings of $N$ and $M$, respectively.
Consequently, $M \cong N$ as braided vector spaces when $\chi
(\alpha)=1,$ i.e. {\rm (i)} holds.

For {\rm (iv)}, assume $\chi (\alpha) =-1.$ If {\rm dim} $\mathfrak
B(\mathcal{O}_{(\alpha, \tau)}^G, \rho) < \infty $, then $q_{M} =
-1$ and $q_{N} = 1$. Consequently, {\rm dim} $\mathfrak
B(\mathcal{O}_{(0, \tau)}^G, \rho) = \infty $. Conversely,  if {\rm
dim} $\mathfrak B(\mathcal{O}_{(0, \tau)}^G, \rho) < \infty $, then
$q_{N} = -1$ and $q_{M} = 1$. Consequently, {\rm dim} $\mathfrak
B(\mathcal{O}_{(\alpha, \tau)}^G, \rho) = \infty $. $\Box$

\begin{Theorem}\label{2.6'}
Let $G = A \rtimes \mathbb S_n$ with $n>2$. Let $\alpha = (1, 1, \cdots, 1), b \in A$ and  $\sigma\in \mathbb
S_n$. If  $(b-\alpha, \sigma)$ can be decomposed   the multiplication of  independent  positive cycles, and
the type of $\sigma$  is not   in the  list below , then {\rm dim}
$\mathfrak B(\mathcal{O}_{(\alpha, \sigma)}^G, \rho) = \infty $ and
{\rm dim}$ \mathfrak B(\mathcal{O}_{(b, \sigma)}^G, \rho') = \infty$ 
for any $\rho \in \widehat {G^{(\alpha, \sigma)}}$, $\rho' \in \widehat {G^{(b, \sigma)}}$.

\renewcommand{\theenumi}{\roman{enumi}}   \renewcommand{\labelenumi}{(\theenumi)}

\begin{enumerate}
\item $(2, 3)$; $(2^3)$;  $(1^m, 2)$;
\item $(3^2)$; $(2^2, 3)$;  $(1^m, 3)$, $(2^4)$; $(1^2, 2^2)$; $(1, 2^2)$;  $(1, p)$; $(p)$ with $p$ prime.

\end{enumerate}
\end{Theorem}

\noindent {\bf Proof.} By Theorem 4.1 of \cite {AFGV}, $\mathcal{O}_{\sigma }^{\mathbb S_n}$ is of  type $D$. 
Obviously, $\mathcal{O}_{\sigma }^{\mathbb S_n}$ is a subrack of $\mathcal{O}_{\sigma }^{G}$, 
hence $\mathcal{O}_{\sigma }^{G}$ is of tpye $D$.
By Lemma \ref {2.6} {\rm (ii)}, $\mathcal{O}_{( \alpha, \sigma)}^{G}=  \alpha \mathcal{O}_{\sigma }^{G}$ is of type $D.$
 It follows from the proof of Theorem 3.6 in \cite  {AFGV} that $\mathfrak B(\mathcal{O}_{(\alpha, \sigma)}^G, \rho) = \infty $.
 Considering Lemma \ref {2.7} (iii), we have the second claim.
 $\Box$

\section { Nichols algebras of reducible {\rm YD} modules over classical Weyl groups} \label {s3}

In this section we present the necessary and sufficient conditions
for Nichols algebras of reducible {\rm YD} modules supported by one
element in $A$ to be finite dimensional. We prove that  if $M$ is a
reducible {\rm YD} module over $ {\mathbb C } G$ supported by $\mathbb S_n$ with
$n \ge 3,$ then ${\rm dim } {\mathfrak B } (M) = \infty$ and if $M$
is a  {\rm YD} module over $ {\mathbb C } G$ supported by $\mathbb A_n$ with $n
\ge 5$,  then  ${\rm dim } {\mathfrak B } (M) =
\infty$.

We say that ${\mathcal O}_{x} $ and ${\mathcal O}_{y}$  are commutative if  $st = ts$ 
for any $s\in {\mathcal O}_{x} $, $t\in {\mathcal O}_{y}$.    
We say that ${\mathcal O}_{x} $ and ${\mathcal O}_{y}$  are squarely commutative if  
\begin {eqnarray}\label {e3.1.1}(st)^2 =( ts)^2
 \end {eqnarray}
 for any $s\in {\mathcal O}_{x} $, $t\in {\mathcal O}_{y}$. This notion is important
according to \cite {HS10} and  will be used in the statements and proofs of Lemmas \ref{3.7} - \ref{3.11'}, Theorem \ref{4'}
and Remark \ref{3.12'}.

\begin {Proposition}\label {3.2}
If $M$ is a reducible {\rm YD} module over ${\mathbb A}_4$ or over
${\mathbb A}_6$,
 then  ${\rm dim } {\mathfrak B } (M) =\infty$.
\end {Proposition}

\noindent {\bf Proof.} Let ${\mathcal F}(G) =: \{ {\mathcal O}_{s}^G  \mid$
$s\in G$,  there exists an irreducible representation $\rho$ of
$G^{s}$ such that $\mathfrak B ({\mathcal O}_{s}, \rho ) < \infty
\}.$

{\rm (i)} Let $G=: {\mathbb A}_4$. It is clear that $s_1=: (1)$,
$s_2=: (1,2)(3,4)$, $s_3=: (1,2,3)$ and  $s_4=: (1,2,4)$ are the
representative system of the conjugacy classes of $\mathbb A_4$.  By Proposition 2.4 of
\cite {AF07}, ${\mathcal O}_{s_2}$ $\notin
{\mathcal F}(G)$. It is clear that $s_2s_3 \not= s_3 s_2$, $s_2s_4
\not= s_4 s_2$, $s_3 (s_2s_3s_2^{-1}) \not= (s_2s_3s_2^{-1})s_3$,
$s_4 (s_2s_4s_2^{-1}) \not= (s_2s_4s_2^{-1})s_4$,  $s_4
(s_4s_3s_4^{-1}) \not= (s_4s_3s_4^{-1})s_4$. Therefore, ${\mathcal
O}_{s_i}$ and ${\mathcal O}_{s_j}$ are not commutative for any $(i,
j) \not = (2, 2)$, $1<i, j \le 4$. It follows immediately that ${\rm
dim } {\mathfrak B } (M) = \infty$ from Theorem 8.2  of \cite {HS10}.

{\rm (ii)}  Let $G=: {\mathbb A}_6$.  If ${\mathcal O}_s \not=
{\mathcal O}_{(12)(3456)}$, then ${\rm dim } {\mathfrak B } (
{\mathcal O}_s, \rho  ) = \infty$ by Theorem 2.7 in \cite {AF07}.
However, ${\mathcal O}_{(12)(3456)}$ and ${\mathcal O}_{(12)(3456)}$
are not commutative since $(12)(34)(12)(3456)(12)(34)= (12)(4356)$
and $(12) (4356) (12) (3456) \not= (12) (3456) (12) (4356)$. It
follows immediately that ${\rm dim } {\mathfrak B } (M) = \infty$
from Theorem 8.2 of \cite  {HS10}. $\Box$

\begin {Theorem}\label {3.3} Let $G= A \rtimes \mathbb S_n$.

{\rm (i)} If $M$ is a reducible {\rm YD} module over $ {\mathbb C } G$ supported
by $\mathbb S_n$ with $n \ge 3,$ then  ${\rm dim } {\mathfrak B }
(M) = \infty$.

{\rm (ii)} If $M$ is a  {\rm YD} module over $ {\mathbb C } G$ supported by
$\mathbb A_n$ with $n \ge 5$, then  ${\rm dim }
{\mathfrak B } (M) = \infty$.

{\rm (iii)} If $M$ is a reducible   {\rm YD} module over $ {\mathbb C } G$
supported by $\mathbb A_n$ with $n \ge 4$,  then ${\rm dim }
{\mathfrak B } (M) = \infty$.
\end {Theorem}

\noindent {\bf Proof.} {\rm (i)} It follows from Corollary 8.4 in \cite {HS10}.
{\rm (ii)} It follows from Theorem 1.1 of \cite {AFGV}.
{\rm (iii)} It follows from Proposition \ref {3.2}. $\Box$

\begin {Theorem}\label {3.4} Let $G= A \rtimes \mathbb S_n$ with $\sigma \in A$ and  $n>2$. Let
$M = M(\mathcal {O}_\sigma ^G, \rho ^{(1)}) \oplus M(\mathcal
{O}_\sigma ^G, \rho ^{(2)}) \oplus \cdots \oplus M(\mathcal
{O}_\sigma ^G, \rho ^{(r)})$ with $\rho ^{(i)} =  \theta _{\chi ,
\mu^{(i)}} = (\chi  \otimes \mu ^{(i)}) \uparrow ^{G^{\sigma_i}}
_{G^{\sigma _i} _{\chi }}$ with $\chi \in \widehat  A$ for $1\le i\le r.$
Then
 {\rm dim} $\mathfrak B(M) < \infty $ if and only if   {\rm dim} $\mathfrak
 B(\mathcal{O}_{\sigma}^G,\rho ^{(i)}) < \infty $ for $1\le i \le r$.

 \end {Theorem}
\noindent {\bf Proof.} The necessity. It is clear. The sufficiency. By Theorem
\ref {2.3}, $f_\sigma$ is odd  and either  $W_\chi = W_\sigma$ or
$|W_\sigma| =n$ or $|W_\chi| =n$. If $f_\sigma$ is odd and
$|W_\sigma| =n$, then it follows from \cite {AS98} that {\rm dim} $\mathfrak B(M) < \infty. $ If
$f_\sigma$ is odd and  $|W_\sigma| <n$, then  $(\mathbb S_n ^\sigma
)_{\chi ^{(i)}}= \mathbb S_n ^\sigma $. Let $P^{(i)}$ and $V^{(i)}$
are representation spaces of $\chi ^{(i)}$ and $\mu ^{(i)}$.
 Let $v_1 ^{(i)}, v_2^{(i)}, \cdots v_{n_i}^{(i)}$ be the basis of $V^{(i)}$,  $0\not=p_0^{(i)}
 \in P^{(i)}$  and $\xi _j ^{(i)}=
  1 \otimes _{k(G^\sigma )_{\chi^{(i)} }} (p_0^{(i)} \otimes v_j^{(i)}) $ for $1\le j \le
  n_i$ with $n_i = {\rm deg} (\mu ^{(i)})$. It is clear that $\{h_i \xi _j^{(\nu)} \mid 1\le i \le m, 1\le
  j \le n_\nu, 1\le \nu \le r \}$ is a basis  of  $M$, where $G = \cup _{i=1}^m h_i G^\sigma
  $ is left coset decomposition of $G^\sigma $ with $h_i \in \mathbb S_n.$
  Let $t_i =: h_i \cdot \sigma = h_i \sigma h_i^{-1}$ and $\gamma _{ij}
  =: h_j ^{-1} h_i \cdot
   \sigma $. Therefore, $t_ih_j = h_j r_{ij}$ and $B(h_i \xi _u ^{(\nu)}\otimes h_j \xi ^{(\nu')}_v )=
   \chi (\gamma _{ij} ) (h_j\xi _v^{(\nu')} \otimes h_i \xi _u^{(\nu)})$. By Lemma \ref {2.2}
   $\chi (\gamma _{ij}) \chi (\gamma _{ji} ) =1$    for any $1\le i, j \le
   m$, $1\le u, v \le  s.$ Consequently, $M$ is of a finite  Cartan type and {\rm dim}
$\mathfrak B(M) < \infty.$ $\Box$

\begin {Theorem}\label {3.6}
Let $G= A \rtimes \mathbb S_n$ and $M = M(\mathcal {O}_{\sigma
_1}^G, \rho ^{(1)}) \oplus M(\mathcal {O}_{\sigma_2} ^G, \rho
^{(2)}) \oplus \cdots \oplus M(\mathcal {O}_{\sigma_r} ^G, \rho
^{(r)})$ with    $\rho ^{(i)} = \theta _{\chi ^{(i)} , \mu^{(i)}} =
(\chi ^{(i)}  \otimes \mu ^{(i)}) \uparrow ^{G^{\sigma_i}}
_{G^{\sigma _i} _{\chi^{(i)} }}$,  $\chi ^{(i)} \in \widehat  A$,
$\sigma _i \in A$ and $n>2$  for $1\le i\le r$. If {\rm dim}
$\mathfrak B(\mathcal{O}_{\sigma_i}^G,\rho^{(i)}) <  \infty$   and
   \begin {eqnarray} \label {e3.6.1}
\mid g(W_{\sigma _i}) \cap W_{\chi ^{(j)}} \mid \equiv
  \mid g^{-1}(W_{\sigma _j}) \cap W_{\chi ^{(i)}} \mid \ \ ( \hbox { \ mod  } \  2)
  \end {eqnarray}
  for any $g \in \mathbb S_n, $ $1\le i, j \le r$, then  {\rm dim} $\mathfrak B(M) < \infty $.
   \end {Theorem}

 \noindent {\bf Proof.} By Theorem \ref {2.3},
$f_{\sigma_i, \chi ^{(i)}}$ is odd and either $W_{\chi ^{(i)}} =
W_{\sigma_i}$ or  $|W_{\sigma_i}| =n$ or  $|W_{\chi ^{(i)}}| =n$ for
$1\le i \le r$.

  If $f_{\sigma_i}$ is odd and $|W_{\sigma_i}|<n$ for any $1\le i\le r$, then
 it follows from Lemma \ref {2.2} that $(\mathbb S_n^{\sigma_i} )_{\chi ^{(i)}}= \mathbb S_n ^{\sigma _i}
 $. Consequently  $\mu ^{(i)} \in \widehat {\mathbb
S_n^{\sigma_i}}$ for $1\le i\le r.$ Let $P^{(i)}$ and $V^{(i)}$ are
representation spaces of $\chi ^{(i)}$ and  $\mu ^{(i)}$.
 Let $v_1 ^{(i)}, v_2^{(i)}, \cdots v_{n_i}^{(i)}$ be the basis of $V^{(i)}$,  $0\not=p_0^{(i)}
 \in P^{(i)}$  and $\xi _j ^{(i)}=
  1 \otimes _{k(G^\sigma )_{\chi^{(i)} }} (p_0^{(i)} \otimes v_j^{(i)}) $ for $1\le j \le
  n_i$ with $n_i = {\rm deg} (\mu ^{(i)})$. It is clear that $\{h_i ^{(\nu)} \xi _j^{(\nu)}
  \mid 1\le i \le m_\nu, 1\le  j \le n_\nu, 1\le \nu \le r \}$ is a basis
  of  $M$, where $G = \cup _{i=1}^{m_\nu} h_i ^{(\nu)} G^{\sigma_\nu}
  $ is left coset decomposition of $G^{\sigma _\nu}$ with $h_i ^{(\nu)} \in \mathbb S_n.$
  Let $t_i ^{(\nu)} =: h_i ^{(\nu)}\cdot \sigma_\nu  = h_i^{(\nu)}
  \sigma _\nu  (h_i^{(\nu)})^{-1}$ and $\gamma _{ij} ^{(u, v)}  =: (h_j ^{(v)})^{-1} h_i^{(u)} \cdot
   \sigma _u$. Therefore, $t_i^{(v)}h_j^{(u)} = h_j ^{(u)}r_{ij}^{(v, u)}$ and
   $B(h_i ^{(u)}\xi _s^{(u)}\otimes h_j^{(v)} \xi _t^{(v)} )=
   \chi ^{(v)}(\gamma _{ij} ^{(u, v)}  ) (h_j^{(v)} \xi _t^{(v)}\otimes h_i ^{(u)}\xi _s^{(u)}
    )$. By condition (\ref {e3.6.1}),
   $\chi ^{(v)}(\gamma _{ij} ^{(u, v)}  ) \chi ^{(u)}(\gamma _{ji} ^{(v, u)}  )=1$
    for any $1\le i\le m _u,$ $1\le  j \le   m_v$, $1\le u, v \le
  r.$ Consequently, $M$ is of a finite  Cartan type and {\rm dim} $\mathfrak B(M) < \infty.$

 Assume that  $f_{\sigma_i}$ is odd for any $1\le i\le r$
and there exists $i_0$ such that  $|W_{\sigma_{i_0}}| =n$. We can
assume that $|W_{\sigma_{i}}| =n$ for $1\le i \le p$ and
$|W_{\sigma_{i}}| <n$ for $p+1\le i \le r$. If $r=p$, then {\rm dim}
$\mathfrak B(M) < \infty$ since it is a central quantum linear space.
 Now $1\le p\le r-1 $. Let $G ^{\sigma _u} =
\cup _{i=1}^{l_u} g_i ^{(u)} (G^{\sigma_u})_{\chi ^{(u)}}
  $ is a left coset decomposition of $(G^{\sigma _u})_{\chi ^{(u)}}$ with $g_i ^{(u)} \in \mathbb
  S_n.$ Set $\xi _{i, j}^{(u)} =: g_i^{(u)} \otimes _{k(G^{\sigma _u})_{\chi ^{(u)}}}
  (p_0^{(u)} \otimes v_j^{(u)})$ for $1\le u\le p$, $1\le i\le l_u$, $1\le j \le n_u$.
  It is clear that  $\{  h_s ^{(u)} \xi _ {i, j}^{( u)}  \mid $ $ 1\le s \le m_u,$ $ 1\le
  u \le p,$  $ 1\le j \le n_u,$ $1\le i \le l_u \}$ $\cup $ $\{  h_i ^{(\nu)} \xi _j^{(\nu)}
  \mid $ $  1\le i \le m_\nu, $ $ 1\le
  j \le n_\nu,$ $ p+1\le \nu \le r \}$ is a basis  of  $M$.  Clearly,
\begin {eqnarray*}
B(h^{(v)}_i \xi _s ^{(v)} \otimes h_j ^{(u)} \xi _{t, z} ^{(u )}) &
=& \chi ^{(u)} ((g_t^{(u)})^{-1} \gamma _{i, j}^{(v, u)}g_t ^{(u)})
(h_j ^{(u)} \xi _{t, z} ^{(u )} \otimes h^{(v)}_i \xi _s^{(v)} )  \ \ \hbox {and }\\
 B(h_j ^{(u)} \xi _{t, z} ^{(u )} \otimes h^{(v)}_i \xi _s
^{(v)} ) & =& \chi ^{(v)} (\gamma _{j, i}^{(u, v)}) (h^{(v)}_i \xi
_s ^{(v)} \otimes h_j ^{(u)} \xi _{t, z} ^{(u )}).
 \end {eqnarray*}

Indeed
 \begin {eqnarray*}
&& f_{(h _i ^{(v)})^{-1}h_j^{(u)} \cdot \sigma _u, \chi ^{(v)}} +
 f_{   (g_t^{(u)})^{-1}(h_j^{(u)})^{-1} h_i^{(v)} \cdot \sigma _v, \chi ^{(u)}} \\
&&~~~~~~\equiv \mid ((g_t^{(u)})^{-1}(h_j^{(u)})^{-1}
 h_i^{(v)})^{-1} (W_{\sigma _u}) \cap W_{\chi ^{(v)}} \mid \\
&& ~~~~~~+  \mid   (g_t^{(u)})^{-1}(h_j^{(u)})^{-1}
 h_i^{(v)} W_{ \sigma _v } \cap W_{\chi ^{(u)}} \mid \
 ( \hbox {since  }  |W_{\sigma _u }| =n)\\
&&~~~~~~\equiv  0 \ \ ( \hbox { mod  } \  2)\ ( \hbox {by condition
} (\ref {e3.6.1}))
 \end {eqnarray*}
and
\begin {eqnarray*} \chi ^{(u)} ((g_t^{(u)})^{-1}
\gamma _{i, j}^{(v, u)}g_t ^{(u)})\chi ^{(v)} (\gamma _{j, i}^{(u,
v)}) &=& (-1) ^{f_{(h _i ^{(v)})^{-1}h_j^{(u)} \cdot \sigma _u, \chi
^{(v)}} + f_{ (g_t^{(u)})^{-1}(h_j^{(u)})^{-1}  h_i^{(v)} \cdot
\sigma _v}} =1
 \end {eqnarray*}
for any $1\le i\le m _v,$ $1\le  j \le   m_u$, $p+1\le v \le r$,$
1\le u \le
  p,$ $1\le s\le n_v$,  $1\le t\le l_u$,  $1\le z \le n_u.$
  Also see
 \begin {eqnarray*}
B(h_j ^{(u)} \xi _{t, z} ^{(u )} \otimes h_{j'} ^{(u')} \xi _{t',
z'} ^{(u' )}) & =& \chi ^{(u')} (\sigma _u) h_{j'} ^{(u')} \xi _{t',
z'} ^{(u' )} \otimes  h_j ^{(u)} \xi _{t, z} ^{(u )} \ \ \hbox {and }\\
\chi ^{(u')} (\sigma _u)\chi ^{(u)} (\sigma _{u'})
 &=& (-1) ^{f_{ \sigma _{u},  \chi ^{(u')}} + f_{ \sigma _{u'}, \chi ^{(u)}}}
 =1\ ( \hbox {by condition } (\ref {e3.6.1})).
 \end {eqnarray*}
   Consequently, considering the proof of Part {\rm (i)}, we have that  $M$ is of a finite
  Cartan type and {\rm dim} $\mathfrak B(M) < \infty.$ $\Box$

\begin {Corollary}\label {3.6''} Let $G= A \rtimes \mathbb S_n$ with $n>2$ and
$M = M(\mathcal {O}_{\sigma _1}^G, \rho ^{(1)}) \oplus M(\mathcal
{O}_{\sigma_2} ^G, \rho ^{(2)}) \oplus \cdots \oplus M(\mathcal
{O}_{\sigma_r} ^G, \rho ^{(r)})$  with  $\sigma _i \in A$ and {\rm
dim} $\mathfrak B(\mathcal{O}_{\sigma_i}^G,\rho^{(i)}) <  \infty$
 for $1\le i \le r$. If  for any $u, v \in \{1, 2, \cdots r\}$  with $u
\not= v,$  one of the following conditions  holds, then  (\ref
{e3.6.1}) holds and {\rm dim} $\mathfrak B(M) < \infty $:

{\rm (i)}   $\sigma _u, \sigma _v \notin Z(G)$
  with $|W_{\chi ^{(u)}}| <n$ and $|W_{\chi ^{(v)}}| <n$.

 {\rm (ii)} $|W_{\chi ^{(u)}}| =n$ and $|W_{\chi ^{(v)}}|= n$.

{\rm (iii)} $|W_{\sigma _u}| =| W_{\sigma _v}| = n$.

{\rm (iv)} $|W_{\chi ^{(u)}}| =n$ and $| W_{\sigma _u}|= n$.

\end {Corollary}
 \noindent {\bf Proof.} Considering Theorem \ref {3.6} we only need show  that (\ref
 {e3.6.1}) holds.

  (a). If {\rm (i)} holds, then $W_{\chi ^{(u)}}= W_{\sigma _u}$ and
 $W_{\chi ^{(v)}}= W_{\sigma _v}$ by Theorem \ref {2.3}. See
\begin {eqnarray*}\mid g(W_{\sigma _u}) \cap W_{\chi ^{(v)}} \mid &=&\mid
g(W_{\sigma _u}) \cap W_{\sigma _v} \mid=\mid g^{-1}(W_{\sigma _v} \cap g ( W_{\sigma _u} ))\mid\\
 &=&  \mid g^{-1}(W_{\sigma _v}) \cap  W_{\chi ^{(u)}} \mid.
  \end {eqnarray*}

(b). Obviously, if {\rm (ii)}  or (iii) holds, then (\ref {e3.6.1})
holds.

 (c). Assume that  {\rm (iv)} holds. Considering (b) and (c), we can assume that  $|W_{\chi
^{(v)}}| <n$ and $|W_{\sigma _v}| < n$. By Theorem \ref {2.3}, $
W_{\chi ^{(v)}}= W_{\sigma _v}$. Thus   (\ref {e3.6.1}) holds.
$\Box$

\begin {Proposition}\label {3.6'''} Let $G= A \rtimes \mathbb S_n$ with $n>2$ and
$M = M(\mathcal {O}_{\sigma _1}^G, \rho ^{(1)}) \oplus M(\mathcal
{O}_{\sigma_2} ^G, \rho ^{(2)}) $  with  $\sigma _2 \in A$ and
$\sigma _2 \not= \sigma _1 =: (1, \cdots, 1)$. If there exists
$g\in \mathbb S_n$ such that (\ref {e3.6.1}) does not hold and some
of the following hold, then  {\rm dim} $\mathfrak B(M) = \infty $:

{\rm (i)}  ${\rm deg \rho ^{(2)}} > 3$.

{\rm (ii)} ${\rm deg \mu ^{(1)}} > 3$.

{\rm (iii)} ${\rm deg \rho ^{(2)}} \ge 2 $ and ${\rm deg \mu
^{(1)}}\ge 2$.

\noindent Here  $\rho ^{(i)} =  \theta _{\chi ^{(i)}, \mu^{(i)}}$
with $\chi ^{(i)}\in \widehat  A$  for $i=1, 2.$
   \end {Proposition}
\noindent {\bf Proof.} If   $f_{\sigma_i}$ is even, then {\rm dim} $\mathfrak
B(M) = \infty $. So we assume that $f_{\sigma_i}$ is odd for $i=1,
2.$  The notations in  the proof of Theorem \ref {3.6} will be used
in the following. Set $u=1$ and $v=2.$ We show this by  the following
several steps.

(a).  If $g\notin G^{\sigma _u} $,  set $g= (g_1 ^{(u)})^{-1}
(h_2^{(u)})^{-1} h_1^{(v)}$ with $h_1^{(v)} = g_1^{(u)} =1.$

(b). If $g\in G^{\sigma _u} $,  set $g= (g_1 ^{(u)})^{-1}
(h_1^{(u)})^{-1} h_1^{(v)}$ with $h_1^{(v)} = g_1^{(u)} =1.$

(c).   We assume that in case (a) $t=1,$ $i = 1$ and $j=2$;  in case
(b) $t=1,$ $i = 1$ and  $j=1$. Clearly,
\begin {eqnarray*}
B(h^{(v)}_i \xi _s ^{(v)} \otimes h_j ^{(u)} \xi _{t, z} ^{(u )}) &
=& \chi ^{(u)} ((g_t^{(u)})^{-1} \gamma _{i, j}^{(v, u)}g_t ^{(u)})
(h_j ^{(u)} \xi _{t, z} ^{(u )}
\otimes h^{(v)}_i \xi _s ^{(v)} ),   \\
 B(h_j ^{(u)} \xi _{t, z} ^{(u )} \otimes h^{(v)}_i \xi _s
^{(v)} ) & =& \chi ^{(v)} (\gamma _{j, i}^{(u, v)}) (h^{(v)}_i \xi
_s ^{(v)} \otimes h_j ^{(u)} \xi _{t, z} ^{(u )})
 \end {eqnarray*} and  \begin {eqnarray*} \chi ^{(u)} ((g_t^{(u)})^{-1}
\gamma _{i, j}^{(v, u)}g_t ^{(u)})\chi ^{(v)} (\gamma _{j, i}^{(u,
v)}) &=& (-1) ^{f_{(h _i ^{(v)})^{-1}h_j^{(u)} \cdot \sigma _u, \chi
^{(v)}} + f_{ (g_t^{(u)})^{-1}(h_j^{(u)})^{-1} h_i^{(v)} \cdot
\sigma _v}} =-1
 \end {eqnarray*} for any  $1\le s\le n_v$,  $1\le z \le n_u.$

 We consider  the Dynkin diagrams of the
braided vector space $ {\mathbb C }  $-span $\{h^{(v)}_i \xi _s ^{(v)},  h_j
^{(u)} \xi _{t, z} ^{(u )} \mid 1\le s\le n_v,  1\le z \le n_u. \}$.
It is clear that there exists a line between every element in
$\{h^{(v)}_i \xi _s ^{(v)} \mid 1\le s\le n_v \}$ and every element
in  $\{  h_j^{(u)} \xi _{t, z} ^{(u )} \mid  1\le z \le n_u \}.$

Obviously, the Dynkin diagrams  of the braided vector space $ {\mathbb C }
$-span $\{h^{(v)}_i \xi _s ^{(v)},  h_j ^{(u)} \xi _{t, z} ^{(u )}
\mid 1\le s\le n_v,  1\le z \le n_u. \}$ are  not of finite Cartan
types in case {\rm (i), (ii), (iii)}. Consequently, {\rm dim}
$\mathfrak B(M) = \infty $ by  Theorem 4 of \cite {He06}. $\Box$

\begin {Lemma}\label {3.7}
Let $G= A \rtimes \mathbb S_n$ with $a \in A, $ $\sigma = (a, \tau )
\in G.$ Then $\mathcal O_a^G$ and $\mathcal O_\sigma ^G$ are
squarely commutative if and only if $a \in Z(G)$ or $\tau ^2=1$
 \end {Lemma}
 \noindent {\bf Proof.} If $\mathcal O_a^G$ and $\mathcal O_\sigma ^G$ are
squarely commutative, then $(\mu \cdot a)\sigma (\mu \cdot a) \sigma
=\sigma (\mu \cdot a)\sigma (\mu \cdot a) $ for any $\mu \in \mathbb
S_n$, which is equivalent to $((\mu \cdot a)a (\tau\cdot (\mu \cdot
a)) (\tau \cdot a), \tau ^2) =(a (\tau\cdot (\mu \cdot a)) (\tau
\cdot a) (\tau ^2\cdot (\mu \cdot a)), \tau ^2) $. It also is
equivalent to
\begin {eqnarray} \label {e3.7.1}
(\mu \cdot a)= \tau ^2\cdot (\mu \cdot a). \end {eqnarray} If $\tau
^2 \not= 1$ and $a \notin Z(G)$, the center of $G$,  then there
exists an independent  cycle  $(i_1, i_2, \cdots, i_r)^{-1}$ of
$\tau ^2$ with $r>1$ and there exist $1\le j _1, j_2 \le n$  such that $a_{j_1} \equiv 0, a_{j_2} \equiv 1$. Thus there exist $\mu \in \mathbb S_n$ such that $ a _{j_1} = (\mu \cdot a)_{i_1} \equiv 0$ and $a _{j_1} =(\mu \cdot a)_{i_2}\equiv
1$. Obviously $(\tau ^2\cdot (\mu \cdot a))_{i_1} = 1$. This
implies $\tau^2 \cdot (\mu \cdot a ) \not= \mu \cdot a. $ Thus
$\mathcal O_a^G$ and $\mathcal O_\sigma ^G$ are not
squarely commutative, which is a contradiction.

Conversely, it is clear  when $a\in Z(G)$. Now assume  that $a\notin
Z(G)$ and $\tau ^2 = 1$. For any $\xi \in \mathbb S_n$, we have
$(\xi \tau \xi ^{-1})^2 =1$. Thus $\mathcal O_a^G$ and $\mathcal
O_\sigma ^G$ are squarely commutative.  $\Box$

\begin {Lemma}\label {3.9'} Let $H = B \rtimes  D$ be a semidirect product of $B$ and $D$, where $B$ is an abelian group. Let $(a, \sigma), (b, \tau)
\in H$ with $a, b \in B$, $\sigma , \tau \in D$. If ${\mathcal
O}_{(a, \sigma)}^H$ and ${\mathcal O}_{(b, \tau)}^H$ are
squarely commutative then ${\mathcal O}_{\sigma}^{D}$ and ${\mathcal
O}_{\tau}^D$ are squarely commutative.
\end {Lemma}

\noindent {\bf Proof.} It is clear that $(a, \sigma ) ^{-1} =( \sigma ^{-1}
\cdot a^{-1}, \sigma^{-1}) $ and \begin {eqnarray}\label {e3.11} (b,
\tau)(a, \sigma) (b, \tau) ^{-1} &=&(b(\tau \cdot a)(\tau \sigma
\tau^{-1} \cdot b^{-1}), \tau \sigma \tau^{-1}).
\end {eqnarray}
For any $x\in {\mathcal O}_{\sigma}^{D}$ and $y\in {\mathcal
O}_{\tau}^D$, by (\ref {e3.11}), there exist $c, d \in B$ such that
$ (c, x) \in {\mathcal O}_{(a, \sigma)}^G$ and $(d, y)\in {\mathcal
O}_{(b, \tau)}^G$. Since $(c, x)(d, y) (c, x)(d, y) = (d, y) (c,
x)(d, y)(c, x)$, we have $xyxy = yxyx$, i.e. ${\mathcal
O}_{\sigma}^{D}$ and ${\mathcal O}_{\tau}^D$ are squarely commutative.
$\Box$

\begin {Lemma}\label {3.10'} Let $1\not= \sigma, 1\not= \tau \in \mathbb S_n $
with $H= \mathbb S_n$ and  $n>2$.  If ${\mathcal O}_\sigma $ and
${\mathcal O}_{\tau}$ are squarely commutative, then one of the
following conditions holds.

{\rm (i)} $n=3$, ${\mathcal O}_{\sigma} = {\mathcal O}_{\tau }=
{\mathcal O}_{(1 2 3)}$ or ${\mathcal O}_{\sigma}={\mathcal O}_{(1
2)}$ and ${\mathcal O}_{\tau}={\mathcal O}_{(1 23)}$.

{\rm (ii)} $n=4$, ${\mathcal O}_{\sigma}= {\mathcal O}_{\tau} =
{\mathcal O}_{(1 2)(3 4)}$ or ${\mathcal O}_{\sigma} = {\mathcal
O}_{(1 2)(3 4)}$ and  $ {\mathcal O}_{\tau} ={\mathcal O}_{(1 23
4)}$ or ${\mathcal O}_{\sigma} = {\mathcal O}_{(12)}$ and $
{\mathcal O}_{\tau} ={\mathcal O}_{(1 2)(3 4)}$.

{\rm (iii)} $n=2k$ with $ {\mathbb C } >2$,  ${\mathcal O}_{\sigma} = {\mathcal
O}_{(12)}$ and $ {\mathcal O}_{\tau} ={\mathcal O}_{(1 2)(3 4)\cdots
(n-1\ n)}$.
\end {Lemma}
\noindent {\bf Proof.} We show this by the following several steps. Assume that
$(1^{\lambda _1}$, $2^{\lambda _2},\cdots , n^{\lambda _n})$ and
$(1^{\lambda _1'}$, $2^{\lambda _2'}, \cdots,  n^{\lambda _n'})$  are the
types of $\sigma $ and  $\tau$, respectively; ${\mathcal O}_\sigma $
and ${\mathcal O}_{\tau}$ are squarely commutative.

{\rm (i)} Let $n = 3$. Obviously,  $\mathcal O_{(1 2 )}$ and
$\mathcal O_{(1 2)}$ are not squarely commutative. Then ${\mathcal
O}_{\sigma} = {\mathcal O}_{\tau }= {\mathcal O}_{(1 2 3)}$ or
${\mathcal O}_{\sigma}={\mathcal O}_{(1 2)}$ and ${\mathcal
O}_{\tau}={\mathcal O}_{(1 23)}$.

{\rm (ii)} Let $n=4$. The types of $\sigma$ and $\tau$ are $(2^2);
(4); (1, 3); (1^2, 2)$, respectively.

(a). ${\mathcal O}_{(1\ 2 ) (3\ 4)}$ and ${\mathcal O}_{(1 \ 2\ 3)}$ are
not squarely commutative since $(12)(34)(123) (12)(34) = (214)$, which
implies $(1\ 2)(3\ 4)(1\ 2\ 3) (1\ 2)(3\ 4)\notin H^{(1\ 2\ 3)}.$

(b). ${\mathcal O}_{(1 234 ) }$ and ${\mathcal O}_{(1234)}$ are not
squarely commutative since \ \ \ $((1234)(4231))^2$ \ \ \ and \ \ \
$((4231)$ $(1234))^2$  map $1$ to $1$ and $2$, respectively.

(c). ${\mathcal O}_{(1 234 ) }$ and ${\mathcal O}_{(123)}$ are not
squarely commutative since $(1234)(123)(1234)$
 maps $1$ to $4$.

 (d). ${\mathcal O}_{(1 234 ) }$ and ${\mathcal O}_{(12)}$ are not
squarely commutative since $(1234)(12)(1234)$
 maps $2$ to $4$.

 (e). ${\mathcal O}_{(1 23 ) }$ and ${\mathcal O}_{(123)}$ are not
squarely commutative since $(134)(123)(134)$
 maps $2$ to $4$.

 (f). ${\mathcal O}_{(1 23 ) }$ and ${\mathcal O}_{(12)}$ are not
squarely commutative since $(14)(123)(14)= (423)$.

(g). ${\mathcal O}_{(1 2 ) }$ and ${\mathcal O}_{(12)}$ are not
squarely commutative since $(14)(12)(14)= (42)$.

{\rm  (iii)} If  $\sigma = (1 2 \cdots r)$  with  $n> 4$,  then
$r=n$ or  ${\mathcal O} _{\sigma} = {\mathcal O} _{(1 2)}$  and
 ${\mathcal O} _{\tau} = {\mathcal O} _{(12) (34) \cdots (n-1\ n)}$ with  $n = 2k >2$.
In fact, obviously,  $H^\sigma $ is a cycle group generated by $(1 2
\cdots r)$.

 (a). If  $\lambda _3' \not=0$ and $n-r
 >1$, set $t= (\cdots, r, r+1, r+2\cdots )t_1$, an independent
 decomposition of $t \in {\mathcal O}_{\tau}$. See $t\sigma t (r) = r+2$, which implies $t\sigma t \notin H^\sigma.$

(b). If  $\lambda _3' \not=0$ and $n-r
 =1$. Set $t= (1,   r, r+1)t_1$, an independent
 decomposition of $t\in {\mathcal O}_{\tau}$. See $t\sigma t (r+1) = t(2)<r+1$ since $r\not= 2$, which implies $t\sigma t \notin H^\sigma.$

 (c). If there exists $j>3$ such that $\lambda _j' \not=0$ with  $n>r>2 $,  set $t= (\cdots, r-2, r-1, r, r+1, \cdots )t_1$, an independent
 decomposition of $t\in {\mathcal O}_{\tau}$. See $t\sigma t (r-2) = r+1$, which implies $t\sigma t \notin H^\sigma.$

(d). If there exists $j>3$ such that $\lambda _j' \not=0$ with
$r=2$, set $t= (1  2  3 \cdots )t_1$, an independent
 decomposition of $t\in {\mathcal O}_{\tau}$. See $t\sigma t (2) = 4$,
  which implies $t\sigma t \notin H^\sigma.$

(e). If the type of $\tau$ is $(1^{\lambda _1'}, 2)$ with $n > r$,
set $t = (r \ r+1) \in {\mathcal O}_{\tau}$.  See $t\sigma t = (1
\cdots r \ r+1) \notin H^\sigma. $

(f). If the type of $\tau$ is $(1^{\lambda _1'}, 2^{\lambda _2'})$ with
$n > r>2$ and $\lambda _2' >1$, set $t= (r \ r+1)(1 2)t_1$, an
independent
 decomposition of $t\in {\mathcal O}_{\tau}$. See $t\sigma t  = (2\ 1 \ \cdots r+1) $
 $ \notin H^\sigma.$

(g). If the type of $\tau$ is $(1^{\lambda _1'}, 2^{\lambda _2'})$ with
$r=2$, $\lambda _2' \ge 1$ and $\lambda _1' \not= 0$, set $t=
(1)(23)t_1$, an independent
 decomposition of $t\in {\mathcal O}_{\tau}$. See $t\sigma t  = (13)$
 $ \notin H^\sigma.$

From now on assume that  both $\sigma$ and $\tau$ are not cycles.

{\rm (iv)} If $n>4$ and $\sigma = (1, 2, \cdots, n)$ is a cycle,
then it is a contradiction.

(a). If  $\lambda _2' \not=0$, set $t= (1, n)(2,3, \cdots )t_1$, an
independent
 decomposition of $t\in {\mathcal O}_{\tau}$.  See $t\sigma t\sigma (1) = t(4)\not=1$ and $\sigma t\sigma t
 (1)=1$.

 (b). If   $\lambda _3' \not=0$, set $t= (1 2 3)t_1$, an independent
 decomposition of $t \in {\mathcal O}_{\tau}$.  See $t\sigma t\sigma  (1) = t(4)$ and $\sigma
 t\sigma t (1) = 2$, which implies that $t\sigma t\notin H^\sigma$
 since $t(4)>3.$

(c). If   $\lambda _4' \not=0$, set $t= (1 2 3 4)t_1$, an
independent
 decomposition of $t \in {\mathcal O}_{\tau}$.  See $t\sigma t\sigma  (1) = 1$ and $\sigma
 t\sigma t (1) = 5$, which implies
 $t\sigma t\notin H^\sigma.$

 (d). If $n>5$ and   there exists $j>4$ such that $\lambda _j' \not=0$, set $t= (1 2 3 4 6 \cdots )t_1$, an independent
 decomposition of $t \in {\mathcal O}_{\tau}$.  See $\sigma t\sigma t (1) = 5$ and $
 t\sigma t\sigma (1) =6$, which implies that $t\sigma t\notin H^\sigma.$

(e). If $n=5$ and   there exists $j>4$ such that $\lambda _j'
\not=0$, set $t= (1  3254 ) \in {\mathcal O}_{\tau}$.  See $\sigma
t\sigma t (1) = 2$ and $
 t\sigma t\sigma (1) =3$, which implies that $t\sigma t\notin H^\sigma.$

{\rm (v)} If $n>4$ and there exists $r>1$ such that $\lambda _r >1$,
then $n = \lambda _rr$.  Let $\sigma = (1, 2, \cdots, r) (r+1,
\cdots, 2r)\cdots ((\lambda _r-1)r+1, \cdots, \lambda _rr)\sigma_1$,
an independent decomposition of $\sigma$. Assume $n > \lambda_r r.$

 (a) If $r>2$ and  $\lambda _j' = 0$ for any $j>2$, set $t = (1, n)(2, 3)t_1$, an
independent decomposition of $t \in \mathcal O_\tau$. See $t\sigma t
= (n, 3, a_3, \cdots, a_ r) \cdots $, which implies that $t\sigma
t\notin H^\sigma.$

(b). If $r=2$ and  $\lambda _j' = 0$ for any $j>2$, set $t = (1,
n)(2, 3)t_1$, an independent decomposition of $t \in \mathcal
O_\tau$. See $t\sigma t = (n, 3)(2, t(4)) \cdots $, which implies
that $t\sigma t\notin H^\sigma$ since $t(4) >2$.

(c). If there exists $j>2$ such that $\lambda _j' \not= 0$,
 set $t = (1, a_1, \cdots, a_p, \lambda _rr, \lambda _rr+1 )(2, 3, \cdots )t_1$, an independent
decomposition of $t \in \mathcal O_\tau$. See $t\sigma t (\lambda
_rr +1)= t\sigma (1)= 3$, which implies that $t\sigma t\notin
H^\sigma.$

{\rm (vi)} If $n>4$ and  $\lambda _r \le 1$ for any $r>1$, then this
is a contradiction. Assume that there exist $r$ and $r'$ such that
$\lambda _{r'} \not= 0$ and $\lambda _{r} \not= 0$ with $2\le r'
<r$. Let $\sigma = (12\cdots r) (r+1 \cdots r+r') \sigma _1$ be an
independent decomposition of $\sigma$.

(a). If $\lambda _i' =0$ for any $i >2$, set $t = (1 \ n) (23) t_1$,
an independent decomposition of $t \in \mathcal O_\tau$. See
$t\sigma t = (n \ 3 \cdots ) \cdots $, which implies that $t\sigma
t\notin H^\sigma.$

(b). If $r> 3$ and there exists $j>2$ such that $\lambda _j'
\not=0$, set $t = (1,  a_1, \cdots, a_p, r, r+1) (2, 3, \cdots )
t_1$, an independent decomposition of $t \in \mathcal O_\tau$. See
$t\sigma t (r+1) =3 $, which implies that $t\sigma t\notin
H^\sigma.$

(c). If $r= 3$ and there exists $j>2$ such that $\lambda _j'
\not=0$, set $t = (123 \cdots ) (34 \cdots ) t_1$, an independent
decomposition of $t \in \mathcal O_\tau$. See $t\sigma t (1) =4 $,
which implies that $t\sigma t\notin H^\sigma.$

{\rm (vii)} If $n>4$ and the types of  $\sigma$ and $\tau$ are  $r
^{\lambda _r}$, then it is a contradiction.  Let $\sigma = (12
\cdots r) (r+1 \cdots 2r) \cdots.$ Set $t= (r+1, 2, \cdots, r) (1,
r+2, \cdots, 2r) t_1$, an independent decomposition of $t \in
{\mathcal O}_{\tau}$. See
 $t\sigma t\sigma (1) =  $ $ \left
\{
\begin
{array} {ll}  5 &\hbox {when } r=3\\
5 &\hbox {when } r>4
\end {array} \right. $
  { and }
 $\sigma t\sigma t(1) = $ $  \left \{
\begin
{array} {ll}
2 &\hbox {when } r=3\\
 r+5 &\hbox {when } r>4
\end {array} \right..$ When $r = 4$, set $t= (4231)t_1$, an independent decomposition of $t \in
{\mathcal O}_{\tau}$. $\sigma t\sigma t (1) = 1$ and $t\sigma t
\sigma (1)= 2$. When $r=2,$ set $t= (14)(25) (36) t_1$, an
independent decomposition of $t \in {\mathcal O}_{\tau}$. See
$\sigma t\sigma t (1) = 5$ and $t\sigma t \sigma (1)= 3$. Then
$\sigma t\sigma t \not= t\sigma t\sigma $. $\Box$

\vskip.1in
In fact, ${\mathcal O}_{ \sigma}$ and ${\mathcal O}_{ \tau}$ in
cases of  Lemma \ref {3.10'} (i) (ii) are squarely commutative .

\begin {Lemma}\label {3.11'}
Let   $G= A \rtimes  \mathbb{S}_n$ with $A\subseteq (\mathbb Z_2)^n$. Then

{\rm (i)}  $n=3$,  ${\mathcal O}_{(a, (12 3))}$ and $ {\mathcal
O}_{(b, (123))}$
 are not squarely commutative.

{\rm (ii)}  $n=3$,    $ {\mathcal O}_{(a, (12))}$  and ${\mathcal
O}_{(b, (12 3))}$
 are not squarely commutative.

{\rm (iii)} $n=4$, $ {\mathcal O}_{(a, (1 2)(3 4))}$ and  ${\mathcal
O}_{(b, (1 23 4))}$ are not squarely commutative.

{\rm (iv)} $n=2k$ with $ k >1$, $ {\mathcal O}_{(a, (1 2))}$ and
${\mathcal O}_{(b, (1 2)(3 4))}$ are not squarely commutative.

{\rm (v)}  $ {\mathcal O}_{(a, \sigma)}$ and ${\mathcal O}_{(b, 1)}$
are  squarely commutative if and only if $\sigma ^2 =1$.

{\rm  (vi)}  If $n=4$, then     ${\mathcal O}_{(a, (12)(34))} $ and
${\mathcal O}_{(b, (12)(34))}$
  are
squarely commutative if and only if   the signs of $(a, (12)(34))$ and
$(b, (12)(34))$ are the same.

\end {Lemma}
\noindent {\bf Proof.} For any $(d, \mu ) \in G= A \rtimes \mathbb{S}_n$, let
$(c, \mu \sigma \mu^{-1})=: (d, \mu) (a, \sigma) (d, \mu)^{-1}$,
i.e. $c = d (\mu \cdot a) (\mu \sigma \mu^{-1} \cdot d)$. It is
clear that $((c, \mu\sigma\mu^{-1}) (b, \tau))^2 = ((b, \tau)(c,
\mu\sigma\mu^{-1}) )^2 $ if and only if
\begin {eqnarray}\label {e3.12}
& & d (\mu \cdot a) (\mu \sigma \mu^{-1} \cdot d)(\mu
\sigma\mu^{-1}\cdot b) (\mu \sigma\mu^{-1}\tau \cdot  d)\nonumber \\
&& (\mu \sigma\mu^{-1}\tau \mu \cdot a) (\mu \sigma\mu^{-1}\tau \mu
\sigma \mu^{-1} \cdot d)  (\mu \sigma\mu^{-1}\tau \mu
\sigma\mu^{-1}\cdot b)   \nonumber \\
&=&  b (\tau \cdot d)
(\tau \mu \cdot a) (\tau \mu \sigma \mu^{-1}\cdot d) (\tau
\mu\sigma\mu^{-1} \cdot b) ( \tau \mu\sigma\mu^{-1}\tau \cdot  d)  \nonumber \\
&& (\tau \mu\sigma\mu^{-1}\tau \mu \cdot a) (\tau \mu\sigma\mu^{-1}
\tau \mu \sigma \mu^{-1} \cdot d),
\end {eqnarray}
which is equivalent to
\begin {eqnarray}\label {e3.13} &&
 d (\mu \sigma \mu^{-1} \cdot d) (\mu \sigma\mu^{-1}\tau \cdot   d
) (\mu \sigma\mu^{-1}\tau \mu \sigma \mu^{-1} \cdot d)  (\tau \cdot d)
 (\tau \mu \sigma \mu^{-1}\cdot d)\nonumber \\
&&( \tau \mu\sigma\mu^{-1}\tau \cdot  d)
(\tau \mu\sigma\mu^{-1} \tau \mu \sigma \mu^{-1} \cdot d)=h
\end {eqnarray}
with $h =: (\mu \cdot a) (\mu \sigma\mu^{-1}\cdot b)
 (\mu \sigma\mu^{-1}\tau \mu \cdot a)  (\mu \sigma\mu^{-1}\tau \mu
\sigma\mu^{-1}\cdot b) b
(\tau \mu \cdot a)  (\tau \mu\sigma\mu^{-1} \cdot b)  (\tau
\mu\sigma\mu^{-1}\tau \mu \cdot a).$

We only need to show that there exists $(d, \mu)\in G$ such that
(\ref {e3.13}) does not hold in the four cases above, respectively.
Let $d= ({d_1}, {d_2}, \cdots, {d_n})$ for any $d\in A.$

{\rm (i)} Let $\sigma = (123)= \tau = \mu$ and $n=3$. (\ref {e3.13})
becomes $d(\sigma \cdot d) =h,$ which  implies $d_1 +d_3 \equiv h_1$
$({\rm mod} \ 2)$. This is a contradiction since $d $ has not this
restriction.

{\rm (ii)} Let $ \tau= (123)$, $ \sigma = (12)=\mu$ and $n=3$. (\ref
{e3.13}) becomes $(\tau^{-1}\cdot d) ((32)\cdot d)(\tau\cdot
d)((13)\cdot d) =h,$ which implies $d_1 +d_2 \equiv h_1$ $({\rm mod}
\ 2)$. This is a contradiction since $d $ has not this restriction.

{\rm (iii)} Let $ \sigma = (12)(34)$, $ \tau = (1234)$, $\mu =
(123)$ and $n=4$. (\ref {e3.13}) becomes $((13)\cdot d) ((4321)
\cdot d) ((1234)\cdot d)((24)\cdot d) =h,$ which implies $d_1
+d_2+d_3+d_4 \equiv h_1$ $({\rm mod} \ 2)$. This is a contradiction
since $d $ has not this restriction.

{\rm (iv)} Let $ \sigma = (12)$, $\lambda = (5 6)(78)\cdots (n-1 \
n)$, $ \tau = (12)(34)\lambda $, $\mu = (123)$ . (\ref {e3.13})
becomes
\begin {eqnarray} \label {e3.11.1}
&&d((23)\cdot d) ((1342)\lambda  \cdot d) ((13)(24)\lambda\cdot
d)((12)(34)\lambda\cdot d) \nonumber\\
&&((1243)\lambda \cdot d) ((14) \cdot d)((14)(23) \cdot d) =h,
\end {eqnarray}
 which implies $0 \equiv h_i \ ({\rm mod} \ 2)$ for $i = 1, \cdots, n$.
 By simple computation,  we have
 $(\mu \cdot a) ((1342) \mu \cdot a) ((12)(34) \mu \cdot a)  ((14) \mu
\cdot a)$ $  ((23)\cdot b) ((13)(24)\cdot b)$
$ b ((1243) \cdot b) = 1,$ which implies $a_3 + a_4 \equiv  0 \ ({\rm
mod} \ 2)$. If $(a, \sigma)$ is a negative cycle, we construct a
negative cycle $(a', \sigma)$ such that $a'_4 \equiv a_3'$ does not
hold as follows: Let $a'_i=0$ when $i\not = 3, 4, $  and $a_4'=1$,
$a_3'=0.$ If $(a, \sigma)$ is a positive cycle, we construct a
positive  cycle $(a', \sigma)$ such that $a'_4 \equiv  a_3'$ does
not hold as follows: Let $a'_i=0$ when $i\not =1, 3, 4, $  and
$a_4'=1$, $a_3'=0$, $a_1' =1$. Since ${\mathcal O}_{(a, \sigma)} =
{\mathcal O}_{(a', \sigma)} $, we obtain a contradiction.

{\rm (v)} Let $\tau =1.$ It is clear that (\ref {e3.13}) becomes
$(\mu \sigma ^2\mu^{-1} \cdot d) d = h$. If $ {\mathcal O}_{(a,
\sigma )}$ and ${\mathcal O}_{(b, 1)}$ are squarely commutative with
$\sigma ^2\not= 1$. then there exists $1\le i \le n$ such that $\mu
\sigma ^{-2}\mu^{-1}(i) = j \not=i$. Thus $d_j +d_i \equiv h_i \
({\rm mod} \ 2)$. This is a contradiction since $d $ has not this
restriction. Conversely, if $   \sigma ^{2}=1 $, then (\ref {e3.13})
holds for any $(d, \mu) \in A \rtimes \mathbb S_n$. Consequently, $
{\mathcal O}_{(a, \sigma )}$ and ${\mathcal O}_{(b, 1)}$ are
squarely commutative.

 {\rm (vi)} It is clear
that ${\mathcal O}_{ (12)(34)} ^{{\mathbb S}_4} $ and ${\mathcal
O}_{ (12)(34)}^{{\mathbb S}_4}$   are
 commutative. (\ref {e3.13}) becomes $1= h$. That is, $(\mu \cdot a)(\mu
\sigma\mu^{-1}\sigma \mu \cdot a)(\sigma \mu \cdot a)(
\mu\sigma\mu^{-1} \mu \cdot a)$ $ (\mu \sigma\mu^{-1}\cdot b) $
$ (\sigma \cdot b)  b   $
$(\sigma \mu\sigma\mu^{-1} \cdot b) $ $  =1$. It is clear that  $\sigma (i)$,
  $\mu \sigma \mu ^{-1}(i)$, $(1)(i)$ and $\mu \sigma \mu ^{-1}\sigma
  (i)$ are different each other for any fixed $i$ with  $1\le i\le 4$ when $\mu \sigma \mu ^{-1}
  \not= \sigma$. Therefore,  (\ref {e3.13}) holds  for any $(d, \mu) \in A\rtimes \mathbb S_n$
  if and  only if  $a_1 +a_2 +a_3 +a_4$ $\equiv$ $b_1 +b_2 +b_3 +b_4$
$({\rm mod } \ 2)$. $\Box$

\begin {Theorem}\label {4'} Let   $G= A \rtimes  \mathbb{S}_n$  with $A\subseteq (\mathbb Z _2)^n$ and  $n>2$. Assume
that there exist different two pairs $(u(C_1), i_1))$ and $(u(C_2),
i_2))$ with $C_1, C_2 $ $\in {\mathcal K}_r(G)$, $i_1 \in I_{C_1}
(r, u)$ and $i_2 \in I_{C_2} (r, u)$ such that $u(C_1)$, $u(C_2)$
$\notin A$. If $ {\rm dim }\mathfrak B (G, r, \overrightarrow{\rho},
u) <\infty $, then the following conditions hold:

 {\rm (i)} $n=4$.

 {\rm (ii)} The type of $u(C)$ is $ (2^{2})$ for any $u(C) \notin A$ and $C\in {\mathcal K}_r(G)$.

{\rm  (iii)} The signs of $u(C)$ and $u(C')$  are the same   for any
$u(C), u(C') \notin A$ and $C, C'\in {\mathcal K}_r(G)$.

\end {Theorem}

\noindent {\bf Proof.} If $ {\rm dim }\mathfrak B (G, r,
\overrightarrow{\rho}, u) <\infty $, then by Theorem
8.6 of \cite{HS10}, ${\mathcal O}_{u(C_1)}^G$ and ${\mathcal O}_{u(C_2)}^G$
are squarely commutative. Let  $ (a, \sigma ) =: u(C_1)$ and $ (b,
\tau ) =: u(C_2)$.  It follows from Lemma \ref {3.9'} that
${\mathcal O}_{\sigma}^{\mathbb S_n}$ and ${\mathcal O}_{\tau
}^{\mathbb S_n}$ are squarely commutative. Considering Lemma \ref
{3.10'}, we have that  one of the following conditions is satisfied

{\rm (i)} $n=3$, ${\mathcal O}_{\sigma}^{\mathbb S_n} = {\mathcal
O}^{\mathbb S_n}_{\tau }= {\mathcal O}^{\mathbb S_n}_{(1 2 3)}$ or
${\mathcal O}^{\mathbb S_n}_{\tau}={\mathcal O}_{(1 2)}$ and
${\mathcal O}^{\mathbb S_n}_{\sigma }={\mathcal O}^{\mathbb S_n}_{(1 23)}$.

{\rm (ii)} $n=4$, ${\mathcal O}^{\mathbb S_n}_{\sigma}= {\mathcal
O}^{\mathbb S_n}_{\tau} = {\mathcal O}^{\mathbb S_n}_{(1 2)(3 4)}$
or ${\mathcal O}^{\mathbb S_n}_{\tau}={\mathcal O}_{(1 2 34)}$ and
${\mathcal O}^{\mathbb S_n}_{\sigma }={\mathcal O}^{\mathbb S_n}_{(1
2)(34)}$. Considering  Lemma \ref {3.11'}, we complete the proof.
$\Box$

\vskip.1in
In other words, we have
\begin {Remark}\label {3.12'} Let $G= A \rtimes \mathbb{S}_n$ with $A\subseteq  (\mathbb Z_2)^n $ and $n>2$.
 Let $M = M({\mathcal O}_{\sigma _1}, \rho ^{(1)})\oplus M({\mathcal
O}_{\sigma _2}, \rho ^{(2)}) \cdots \oplus M({\mathcal
O}_{\sigma_m}, \rho ^{(m)})$ be a reducible {\rm YD} module over
$ {\mathbb C } G$. Assume that  there exist $i\not= j$  such that $\sigma _i$,
$\sigma _j $ $\notin A$. If  $ {\rm dim }\mathfrak B (M)<\infty$,
then  $n = 4$,  the type of $\sigma _p$ is $(2^2)$ and the sign of
$\sigma _p$ is stable  when $\sigma _p \notin A.$

\end {Remark}

Note that we use Lemma \ref {3.9'} in proof of this Theorem \ref {4'}, i.e. we prove  this theorem
by means of squarely commutative property.  Thus we
  need Lemma \ref {3.10'} in stead of Corollary 8.4 in \cite {HS10}.

\begin {Theorem}\label {3.8} Let $G= A \rtimes \mathbb S_n$ with
$n>2$ and $M = M(\mathcal {O}_{\sigma _1}^G, \rho ^{(1)}) \oplus
M(\mathcal {O}_{\sigma_2} ^G, \rho ^{(2)}) \oplus \cdots \oplus
M(\mathcal {O}_{\sigma_m} ^G, \rho ^{(m)})$ be a reducible {\rm YD}
module.
 If {\rm dim} $\mathfrak B(M) < \infty $ then one of  the following
 holds:

{\rm (i)} $m_1 \ge m-1$.

{\rm (ii)}   $m_2 = m-1$ and $ \sigma _i = (a, \tau )$ with
$\tau^2=1$, when  $\sigma_i \notin A$.

{\rm (iii)} $m_3 \ge  2$,  $n=4$,  $ \mathcal O_ {\sigma _i}=
\mathcal O_ {(a, (12)(34))}$; the signs of $\sigma _i $ and $\sigma
_j$ are the same; $\sigma_ i$ has a negative cycle,  when  $\sigma _i, \sigma _j \notin A$.

Here $m_1 =: \mid \{ i  \mid \sigma _i =(1, 1, \cdots, 1),
1\le i \le m \}\mid$,  $m_2 =: \mid \{ i  \mid  \sigma_i \in A, 1\le
i \le m \}\mid $,  $m_3 =: \mid \{ i \mid \sigma _i \notin A, 1\le i \le m \}\mid$.

 \end {Theorem}
 \noindent {\bf Proof.} Let $\alpha =: (1, 1, \cdots, 1)$.
If  $m_2 = m-1$, then  $\tau ^2=1$ by Lemma \ref {3.7}. If  $m_2 <
m-1$, then case {\rm (iii)} holds  by Remark \ref  {3.12'} and Theorem \ref {2.7'}. $\Box$

\vskip.1in
 A central quantum linear
space is a finite dimensional Nichols algebra, which was introduced
in \cite {AS98,AS02}. ${\rm RSR}(G, r, \overrightarrow
{\rho}, u)$ is said to be  a {\it  central
 quantum linear type}  if it is  quantum  symmetric
  and  of the non-essentially infinite type with
$C \subseteq Z(G)$ for any $C\in {\mathcal K}_r(G)$. In this case,
$\mathfrak {B} (G, r, \overrightarrow \rho, u)$ is called a {\it
central  quantum linear space} over $G$.

\begin {Theorem}\label {3'}
$\mathfrak {B}  (G, r, \overrightarrow {\rho}, u)$   is a central
quantum linear space over  classical Weyl group $G$ if and only if
  $C= \{  (1, \cdots, 1) \}\subseteq G$, $r = r_CC$, $\rho_C^{(i)}= \theta _{\chi _C^{(i)}, \mu_C^{(i)}} =:
  (\chi _C^{(i)} \otimes \mu^{(i)}_C) \uparrow _{G_{\chi _C^{(i)}} ^{u(C)}}^{G^{u(C)}}\in  \widehat
  {G^{u(C)}}$
  with $\chi _C^{(i)}$ $ \in \{  \chi _2^{\delta ^{(i)}_1} \otimes \chi _2^{\delta ^{(i)}_2}
\otimes \cdots \otimes \chi _2^{\delta ^{(i)}_n}$ $\mid $  $\delta
^{(i)}_1 +\delta ^{(i)}_2 + \cdots + \delta ^{(i)}_n$ is odd \} for
any $i \in I_C(r, u)$.
\end {Theorem}
\noindent {\bf Proof.} It is clear $\theta _{\chi _C^{(i)}, \rho_C^{(i)}} (
(1, \cdots, 1)) = \chi_2 (1)^{\delta ^{(i)}_1 +\delta ^{(i)}_2
+ \cdots + \delta ^{(i)}_n }\ {\rm id} $.  Applying \cite {AS98}, we
complete the proof. $\Box$

\vskip.1in
In other words we have

\begin {Remark}\label {3.8'} Let $G= A \rtimes \mathbb{S}_n$.
 Assume that  $\alpha =: (1, 1,\cdots, 1) \in G$ and
$M = M({\mathcal O}_\alpha , \rho  ^{(1)})\oplus M({\mathcal
O}_\alpha, \rho ^{(2)}) \oplus \cdots \oplus M({\mathcal O}_\alpha,
\rho ^{(m)})$ is a reducible
 {\rm YD} module over $ {\mathbb C } G$ with $ \rho  ^{(i)} =
\theta _{\chi ^{(i)}, \mu ^{(i)}} =:
  (\chi ^{(i)} \otimes \mu ^{(i)}) \uparrow _{G_{\chi{(\nu_i)}}^a }^{G^a}\in  \widehat
  {G ^{a}}$ for $i=1, 2, \cdots, m$. Then
$\mathfrak B (M)$ is  finite dimensional if and only if
$\chi^{(i)}(\alpha) = -1$ for $1\le i \le m$.
  \end {Remark}
In fact, if there is $1\le i \le m $, such that $\chi ^{(i)}
(\alpha) =1$, then  ${\rm dim} \mathfrak B (M)= \infty$ by
Proposition \ref {2.6}.

\section {Relationship between $W(B_n)$ and $W(D_n)$} \label {s4}

In this section we establish the relationship between Nichols algebras over
$W(B_n)$ and $W(D_n)$.
\begin{Lemma}\label {4.1} (Remark 2.10 in \cite {AF07})
 Let $H$ be a subgroup of group $D$ with  index two.  Assume that  $\sigma \in
 H$ with ${\mathcal O}_\sigma ^D={\mathcal O}_\sigma ^H$ and
$\eta \in \widehat {D^\sigma }$. Then there exists $\rho \in
\widehat {H^\sigma }$ such that $M(\mathcal{O}_\sigma^D, \eta) \cong
M(\mathcal{O}_\sigma^H, \rho)$ or $M(\mathcal{O}_\sigma^D, \eta)
\cong  M(\mathcal{O}_\sigma^H, \rho) \oplus M(\mathcal{O}_\sigma^H,
\overline{\rho}) $ as braided vector spaces.

\end {Lemma}
\noindent {\bf Proof.} It is easy to see that
\begin{eqnarray} \label {e4.1}
{\eta}'(g): =
\begin{cases} \eta(g), & \hbox {  if  }   g \in {H^\sigma}, \\ -\eta(g) & \hbox { if  }  g \in  D^\sigma
\setminus H ^\sigma ,
\end{cases}
\end{eqnarray}
defines a new representation of $\mathrm{D}$.

We have two cases (see Subsection 2.5 of \cite {AF07}):

{\rm (i)} $\eta \ncong \eta'$. If $\rho : =\eta \mid _{H^\sigma}$,
then $\rho \in \widehat{H^\sigma}$, $\rho \cong \overline{\rho}$ and
$ \rho \uparrow  _{H^\sigma}^{D^\sigma}\cong \eta \oplus \eta'$.

{\rm (ii)} $\eta \cong \eta'$. We have that $ \eta \mid _{H^\sigma}
\cong \rho \oplus \overline{\rho}$ and $ \rho \uparrow
_{H^\sigma}^{D^\sigma} \cong \eta \cong  \overline{\rho} \uparrow
_{H^\sigma}^{D^\sigma}$.

$M(\mathcal{O}_\sigma^D, \eta) \cong M(\mathcal{O}_\sigma^H, \rho)$,
\ for the case {\rm (i)}, as braided vector spaces.

$M(\mathcal{O}_\sigma^D, \eta) \cong M(\mathcal{O}_\sigma^H, \rho)
\oplus M(\mathcal{O}_\sigma^H, \overline{\rho}) $, \ for the case
{\rm (ii)}, as braided vector spaces.

 By Proposition 2.3 of \cite {Se77}, $\rho$ is irreducible. $\Box$

\begin {Proposition}\label {4.2} Let $G =: A \rtimes \mathbb S_n = W(B_n)$
and $H =: W(D_n)$ with $n>3$. Assume that  $\sigma = (a, \xi ) \in
W(D_n)$ such that $G^\sigma = A^\sigma \rtimes \mathbb S_n^\sigma $.

{\rm (i)} If there exists a negative sign cycle in $(a, \xi)$  or
there exists  a cycle in $(a, \xi)$ such that its length is odd,
then $\rho \in \widehat {H^\sigma}$ and {\rm dim}$\mathfrak B
(\mathcal{O}_\sigma^G, \eta ) = {\rm dim} \mathfrak
B(\mathcal{O}_\sigma^H, \rho)$ for any $\eta \in \widehat
{G^\sigma}$ with $\rho = \eta \mid _{H^\sigma}$.

{\rm (ii)} If every sign cycle of $(a, \xi)$ is positive and the
length of every cycle is even, then {\rm dim}$\mathfrak B
(\mathcal{O}_\sigma^G, \rho ) = \infty$ when  ${\rm dim} \mathfrak
B(\mathcal{O}_\sigma^H, \rho) = \infty$ with  $\rho \in \widehat {G^\sigma}$.

{\rm (iii)} If $n$ is odd,
then $\rho \in \widehat {H^\sigma}$ and {\rm dim}$\mathfrak B
(\mathcal{O}_\sigma^G, \eta ) = {\rm dim} \mathfrak
B(\mathcal{O}_\sigma^H, \rho)$ for any $\eta \in \widehat
{G^\sigma}$ with $\rho = \eta \mid _{H^\sigma}$.

\end{Proposition}

\noindent {\bf Proof.} {\rm (i)} By Lemma \ref {6.1'} in the Appendix, $\mathcal{O}_\sigma^G=
\mathcal{O}_\sigma^H$. Let $S= \{ a \in A \mid  a= ({a_1}, {a_2}, \cdots, {a_n})\}$. It is clear that $W(D_n) = S\rtimes
\mathbb S_n$.  Now we  show that $A^\sigma \nsubseteq S$. In
fact, if $A^\sigma \subseteq S$, then $G^\sigma = H^\sigma$, which
is a contradiction. Let $\beta \in A^\sigma \setminus S.$

 We first show that $\eta \ncong \eta'$ for
any $\eta \in \widehat {G^\sigma}$,  where $\eta'$ is defined in
(\ref {e4.1}). In fact, assume $\eta \cong \eta'.$ Then $\chi _\eta (x) = 0$ for any $x\in
G^\sigma \setminus H^\sigma$, where $\chi _\eta $ denotes the
character of $\eta$. On the other hand, $\chi _\eta (x) = [\mathbb
S_n^\sigma , (\mathbb S_n^\sigma) _\chi] \chi _{\chi \otimes \mu}
(x)$ since $\eta = (\chi \otimes \mu ) \uparrow ^{G^\sigma}
_{(G^\sigma)_\chi}$ by Definition 1.12.2 in \cite {Sa01}. Consequently,
$\chi _{\chi \otimes \mu} (x)=0$. For any $\tau \in  (\mathbb
S_n^\sigma)_\chi$, we have that $(\beta, \tau)\in G^\sigma \setminus
H^\sigma $. See $\chi _{\chi \otimes \mu} (\beta, \tau) = \chi
(\beta) \chi _\mu (\tau)=0$. Therefore, $\chi _\mu (\tau)= 0$, which
is a contradiction since $\chi _\mu (1 ) \not=0$. By (i) in the
proof of Lemma \ref {4.1}, $M(\mathcal{O}_\sigma^G, \eta) \cong
M(\mathcal{O}_\sigma^H, \rho)$. Consequently, {\rm dim}$\mathfrak B
(\mathcal{O}_\sigma^G, \eta ) = {\rm dim} \mathfrak
B(\mathcal{O}_\sigma^H, \rho)$.

{\rm (ii)} By Lemma \ref {6.1'} in the Appendix, $G^\sigma =
H^\sigma.$ It follows from Lemma 2.2 of \cite  {AFGV}.

{\rm (iii)} It follows from {\rm (i)}.
$\Box$


\section{Appendix} 

In this Appendix the conjugacy  classes of Weyl groups of types $B_n$
and $D_n$ are found. They were obtained in \cite {Ca72} (Pages
25-26). Here for completeness we give  proof of the results
by using semi-direct products, while  in \cite {Ca72} admissible diagrams were used.

\begin{Lemma}
Let $G= A \rtimes \mathbb S_n$ with $A = (\mathbb Z_2 )^n$ and $n>2$.
Then $W(B_n) \cong A \rtimes \mathbb S_n.$
\end {Lemma}

  Let $G= A \rtimes \mathbb S_n$. $(a, \sigma)\in G$ is called a
sign cycle if $\sigma = (i_1, i_2, \cdots, i_r)$  is cycle and $a =
({a_1}, \cdots, {a_n})$ with $a_i=0$ for $i \notin \{i_1,
i_2, \cdots, i_r\}$.  A sign cycle $(a, \sigma)$ is called positive
( or negative ) if $\sum _{i=1}^n a_i $ is even (or odd). $(a,
\sigma ) = (a^{(1)}, \sigma _1)(a^{(2)}, \sigma _2)\cdots (a^{(r)},
\sigma _r)$ is called an independent sign cycle decomposition of
$(a, \sigma)$ if $\sigma = \sigma _1 \sigma _2 \cdots \sigma _r$ is
an independent cycle decomposition of $\sigma$ in $\mathbb S_n$ and
$(a^{(i)}, \sigma _i)$ is a sign cycle for $1\le i\le r$. $(a, \sigma)$ is called to be positive (negative ) 
if the number of negative cycles is even (odd).
Sign cycle decompositions exist and are unique up to rearrangement.

\begin {Lemma}\label {6.0}  $W(D_n) = \{ (a, \tau) \in W(B_n)\mid$ $
\ \hbox { the number of negative cycles  in }$ $ (a, \tau) \ $ $
\hbox {is even } \}= \{ (a, \tau)\in W(B_n) \mid$ $ \ \hbox { the
sign of  }$ $ (a, \tau)  $ $ \hbox {is positive } \}$.
\end {Lemma}
\noindent {\bf Proof.} $e_i - e_j,$ $ -(e_i+e_j),(e_i+e_j)$, for $i\not= j, i,
j =1, 2, \cdots, n $,  are roots of  $D_n$ with orthonormal
 set $\{e_i \mid 1\le i \le n\}$. Furthermore, $ e_1-e_2,$ $ e_2 - e_3,$ $ \cdots, $
$ e_{n-1} - e_n$, $e_{n-1}+e_n$ is a prime root system of $D_n$. Let
$r_i =: (i, i+1)$, $\tau _i =: \sigma _{e_i}$  and $\bar r_i =:
\sigma _{e_i-e_{i+1}}$. If $a = ({a_1}, \cdots, {a_n})$ with
$a_j =1$ when $j=i$  and $a_j =0$ otherwise, then $a$ is denoted by
$\widetilde{ \tau_i}$. It is clear $\sigma _{e_{n-1}+e_{n}} (e_s) =
\tau _{n-1} \tau _n e_{r_{n-1} (s)}$ and $\sigma _{e_{n-1}+e_{n}} =
\tau _{n-1}\tau _n \bar r_{n-1} $. It is enough to show the numbers
of negative cycles in $r_i (a, \mu) $ and  $\widetilde{\tau
_{n-1}}\widetilde{\tau _n} r_{n-1} (a, \mu) $   are  even for $1\le
i \le n-1$, respectively,  when one of $(a, \mu)$ is even  since
$W(D_n)$ is generated by $r_i$, $\widetilde{\tau
_{n-1}}\widetilde{\tau _n }r_{n-1}$. We only need show that the
number of negative cycles in $r_i (a, \mu) $ is even for $1\le i \le n-1$.

Assume $a= ({a_1}, \cdots, {a_n})$. By (\ref {e6.1.2}),
$r_i (a, \mu) = ({a_{r_i(1)}}, \cdots, {a_{r_i(n)}})r_i\mu$

{\rm (i)} Assume that there exists a cycle $(i_1, i_2, \cdots, i_s)$
of $\mu$ such that $i = i_1 $ and $i_t = i+1$. Then $r_i(i_1, i_2,
\cdots, i_s) =  (i_1, i_2, \cdots, i_{t-1})(i_{t}, i_{t+1}, \cdots,
i_{s})$. If $(i_1, i_2, \cdots, i_s)$  is a negative cycle in $(a,
\mu)$, then the signs of $(i_1, i_2, \cdots, i_{t-1})$ and $(i_{t},
i_{t+1}, \cdots, i_{s})$ in $r_i (a, \mu) $ are not the same.
 If $(i_1, i_2, \cdots, i_s)$  is a
positive  cycle of $(a, \mu)$, then the signs of $(i_1, i_2, \cdots,
i_{t-1})$ and $(i_{t}, i_{t+1}, \cdots, i_{s})$ in $r_i (a, \mu) $ are  the same.

{\rm (ii)} Assume that there exist two different independent  cycles
$(i_1, i_2, \cdots, i_s)$  and $(j_1, j_2, $ $ \cdots, j_t)$ of
$\mu$ such that $i = i_1 $ and $j_1 = i+1$. Then $r_{i}(i_1, i_2,
\cdots, i_s)$ $(j_1, j_2, $ $ \cdots, j_t)$ $=   (i_1, i_2, $
$\cdots, i_s, j_1, j_2, $ $ \cdots, j_t)$. If the signs of $(i_1,
i_2, \cdots, i_s)$ and $(j_1, j_2, \cdots, j_t)$ in $(a, \mu)$ are
the same, then the sign of $(i_1, i_2, \cdots, i_s, j_1, j_2,
\cdots, j_t)$ is positive. If the signs of $(i_1, i_2, \cdots, i_s)$
and $(j_1, j_2, \cdots, j_t)$ in $(a, \mu)$ are not the same,  then
the sign of $(i_1, i_2, \cdots, i_s, j_1, j_2, \cdots, j_t)$ is negative.

(iii) Consequently, the number of negative cycles in $r_i (a, \mu) $
is even if and only if the number of negative cycles in $(a, \mu)$
is even. $\Box$

\begin {Lemma}\label {6.1} Let $G= A \rtimes \mathbb S_n$. Assume that $(a, \sigma)$ and $(a', \sigma')$ are
two sign cycles. Then $(a, \sigma)$ and $(a', \sigma')$ are
conjugate in $G$ if and only  if the lengths of $\sigma$ and $
\sigma'$ are the same and $\sum _{i=1}^n a_i \equiv \sum _{i=1}^n
a'_i  \ \ (\hbox {mod} \ \ 2)$.
\end {Lemma}
\noindent {\bf Proof.} If  $(a, \sigma)$ and $(a', \sigma')$ are conjugate in
$G$, then there exists $(b,\tau ) \in G$ such that $(b, \tau) (a,
\sigma) (b, \tau)^{-1}  = (a', \sigma')$. By simple computation we
have \begin {eqnarray} \label {e6.1.1} a' = b(\tau \cdot a)(\sigma
'\cdot b)  \ \ \hbox { and } \ \ \tau \sigma \tau ^{-1} =
\sigma'.\end {eqnarray} Thus $ \sigma$ and $\sigma'$ are conjugate
in $\mathbb S_n$, which implies  that the lengths of $\sigma$ and $
\sigma'$ are the same. Furthermore, $\sum _{i=1}^n a_i \equiv \sum
_{i=1}^n a'_i \ \ (\hbox {mod} \ \  2)$.

Conversely, assume  the lengths of $\sigma$ and $ \sigma'$ are the
same and $\sum _{i=1}^n a_i \equiv \sum _{i=1}^n a'_i  \ \ (\hbox
{mod} \ \ 2) $ with  $\sigma = (i_1, i_2, \cdots, i_r)$ and $\sigma'
= (i_1', i_2', \cdots, i_r')$. Let $\tau \in \mathbb S_n$ such that
$\tau (i _j) = i'_j$ for $1\le j \le r$. It is clear that $a' = \tau
\sigma \tau ^{-1}.$ It is enough to
 find $b \in A$ such that (\ref {e6.1.1}) holds. By (\ref
{e6.1.1}), we have $a_{i_j'}' \equiv b_{i_j'} + a_{\tau
^{-1}(i_j')}+ b _{\sigma '{}^{-1}(i_j')} \ \ (\hbox {mod} \ \ 2)$
for $1\le j \le r$. This implies that 
$$\left \{ \begin {array} {lll}
a_{i_1'}' &\equiv & b_{i_1'} + a_{i_1}+ b _{i_r'} \ \ \ (\hbox {mod} \ \ 2)\\
a_{i_2'}' &\equiv & b_{i_2'} + a_{i_2}+ b _{i_1'} \ \ \ (\hbox {mod} \ \ 2) \\
\cdots & & \cdots   \\
a_{i_r'}' &\equiv & b_{i_r'} + a_{i_r}+ b _{i_{r-1}'} \ \   (\hbox
{mod}\ \  2). \end {array} \right.$$
That is, 
$$ \left ( \begin
{array} {llllll}
1 & 0 & 0& \cdots & 0 & 1 \\
1&1 & 0&\cdots & 0 & 0\\
0 &1 & 1&\cdots & 0 & 0\\
\cdot &\cdot &\cdot &\cdots & \cdot & \cdot\\
0&0 & 0&\cdots & 1 & 1
 \end {array} \right) \left ( \begin {array} {l}
b_{i_1'}\\
b_{i_2'}\\
b_{i_3'}\\
\cdot \\
b_{i_r'}
 \end {array} \right) \equiv  \left ( \begin {array} {l}
a_{i_1'}' + a_{i_1}\\
a_{i_2'}' + a_{i_2}\\
a_{i_3'}' + a_{i_3}\\
\cdot \\
a_{i_r'}' + a_{i_r}
 \end {array} \right)  (\hbox {mod} \ \ 2).$$  
The solutions are
\begin {eqnarray} \label {e6.1.2}\left ( \begin {array} {l}
b_{i_1'}\\
b_{i_2'}\\
b_{i_3'}\\
\cdot \\
b_{i_{r-1}'}\\
b_{i_r'}
 \end {array} \right) \equiv  \left ( \begin {array} {l}
a_{i_1'}' + a_{i_1} + b_{i_r'}\\
a_{i_1'}' + a_{i_1} +a_{i_2'}' + a_{i_2}+  b_{i_r'}\\
 a_{i_1'}' + a_{i_1}+ a_{i_2'}' + a_{i_2}+ a_{i_3'}' + a_{i_3} +b_{i_r'}\\
\cdots \\
 a_{i_1'}' + a_{i_1}+ a_{i_2'}' + a_{i_2}+ \cdots + a_{i_{r-1}'}' +
 a_{i_{r-1}}+b_{i_r'}\\
 b_{i'_r}
 \end {array} \right) (\hbox {mod} \ \ 2), \end {eqnarray} where $b_{i_r'} =
 1$ or $0.$
$\Box$

\begin {Theorem}\label {6.2} Let $G= A \rtimes \mathbb S_n$. Assume
$(a, \sigma ) = (a^{(1)}, \sigma _1)(a^{(2)}, \sigma _2)\cdots
(a^{(r)}, \sigma _r)$ and $(a', \sigma ') = (a'{}^{(1)}, \sigma'
_1)(a'{}^{(2)}, \sigma' _2)\cdots (a'{}^{(r')}, \sigma' _{r'})$ are
independent sign cycle decomposition of $(a, \sigma)$ and $(a',
\sigma')$, respectively. Then  $(a, \sigma)$ and $(a', \sigma')$ are
conjugate in $G$ if and only  if $(a^{(i)}, \sigma_i)$ and
$(a'{}^{(i)}, \sigma'_i)$ are conjugate in $G$  for $1\le i \le r$
with  $r= r'$  after rearranging.
\end {Theorem}
\noindent {\bf Proof.} The sufficiency. Assume that $(a^{(i)}, \sigma_i)$ and
$(a'{}^{(i)}, \sigma'_i)$ are conjugate in $G$  for $1\le i \le r$
and $r=r'$. It follows from the proof of Lemma \ref {6.1} that there
exist $\tau \in \mathbb S_n$ and $b\in A$ such that $\tau \sigma _i
\tau^{-1} = \sigma _i'$ and $b (\tau \cdot a^{(i)}) (\tau \sigma _i
\tau ^{-1} \cdot b) = a'{}^{(i)}$ for $1\le i \le r$. Consequently,
$(b, \tau)(a, \sigma )(b, \tau^{-1})$ $ = (a' , \sigma')$.

The necessity. There exists $(b, \tau)\in G$ such that $(b, \tau)(a, \sigma )(b, \tau^{-1})$ $ =
(a' , \sigma')$. Since $\tau \sigma \tau ^{-1} = \sigma'$, we have
$r=r'$. It is enough to show that $ (a^{(i)}, \sigma _i)$ and $(a'{}^{(i)}, \sigma _i')$
 are conjugate for $1\le i \le r$ after rearranging. It follows from Lemma \ref {6.1}.
$\Box$

\vskip.1in
It is clear that $W(B_n) = W(D_n) \cup \beta W(D_n)$ is a left coset
decomposition of $W(D_n)$ in $W(B_n)$, where  $\beta = (1, 0, \cdots, 0)\in A$.

\begin {Lemma}\label {6.1'}
Let $(a, \mu) \in W(D_n)$ with $n\ge  4.$

{\rm (i)} If there exists a negative sign cycle in $(a, \mu)$  or
there exists  a cycle in $(a, \mu)$ such that its length is odd,
then $W(B_n) ^{(a, \mu)} \not= W(D_n)^{(a, \mu)}$ and $\mathcal O_{(a, \mu)} ^{W(B_n)} = \mathcal O_{(a, \mu)}^{W(D_n)}$.

{\rm (ii)} If every sign cycle of $(a, \mu)$ is positive and the
length of every cycle is even, then  $W(B_n) ^{(a, \mu)} =W(D_n)^{(a, \mu)}$.
\end {Lemma}
\noindent {\bf Proof.} Let the type of $\mu$ be $(1^ {\lambda _1}, 2^{\lambda_2}, \cdots,  n ^{\lambda _n})$.

{\rm (i)} Assume that  there exists a sign cycle $(c, \xi)$ in $(a,\mu)$
 with $\xi = (i_1, i_2, \cdots, i_r)$ and odd $r$.  Set $\tau =1$ and
$b= ({b_1},{b_2}, \cdots, {b_n} ) $ with $b_i = 1 $
when $i = i_1, i_2, \cdots, i_r$, otherwise  $b_i =0$. Obviously,
$(b, \tau) \in W(B_n) ^{(a, \mu)}$ and $(b, \tau) \notin W(D_n) ^{(a,
\mu)}$. Therefore, $W(B_n) ^{(a, \mu)} \not= W(D_n) ^{(a, \mu)}$ and
$\mathcal O_{(a, \mu)} ^{W(B_n)} = \mathcal O_{(a, \mu)} ^{W(D_n)}$.

Assume that   there exists  a cycle $(c, \xi)$ in $(a, \mu)$ such
that $(c, \xi)$ is negative with $\xi = (i_1, i_2, \cdots, i_r)$.
Therefore, $W(B_n) ^{(a, \mu)} \not= W(D_n) ^{(a, \mu)}$ and
$\mathcal O_{(a, \mu)} ^{W(B_n)} = \mathcal O_{(a, \mu)} ^{W(D_n)}$.

We prove {\rm (ii)} by the following several steps.

(a).  Assume  $\mu = (i_1i_2\cdots i_n)$ with $n=2l$. If $(b, \tau)
\in W(B_n)^{(a, \mu)}$ with $i_j' =: \tau (i_j) $ for $1\le j \le n$, then
\begin {eqnarray} \label {e6.1'.1}\sum _{j=1}^nb_{i_j'} = \sum
_{j=1} ^l a_{i_{2j -1}'} + \sum _{j=1} ^l a_{i_{2j -1}} \end
{eqnarray} by (\ref {e6.1.2}) and $\tau = \mu ^k$ for some $ {\mathbb C } .$ It
is clear that $\sum _{j=1} ^l a_{i_{2j -1}'} = \sum _{j=1} ^l
a_{i_{2j -1}}$ when $\tau = \mu^{2k}$ and $\sum _{j=1} ^l a_{i_{2j
-1}'} = \sum _{j=1} ^l a_{i_{2j}}$ when $\tau = \mu^{2k-1}$.
Consequently, $(b, \tau) \in W(D_n)^{(a, \mu)}$, i.e. $W(B_n) ^{(a,\mu)} =W(D_n) ^{(a, \mu)}$.

(b). Assume that the type of $\mu$ is $(r)^{\lambda _{r}}$ with even
$r$ and $\mu = (1, 2, \cdots, r)(r+1, $ $ \cdots, 2r) $ $\cdots
((\lambda _r-1)+1, (\lambda _r-1)r+2,  $ $ \cdots, \lambda _r r)$.
Then $\mathbb S_n ^{\mu}$ is generated by $A_{k,r}$ and $B_{h, r}$
for $1\le k \le \lambda _r$, $1\le h\le \lambda _r-1,$ where $A_{k,
r} = ((k-1)r +1, (k-1)r +2, \cdots, kr )$ and $B_{h, r} = ((h-1)r
+1, hr+1)((h-1)r +2, hr+2)\cdots ((h-1)r +r, hr+r)$.

(c). If the sign cycle $(b, \tau ) \in W(B_n)^{(a, \mu)}$ with $\tau =
A_{k, r}$, then $b(\tau \cdot a) (\mu \cdot b) = a$. This implies that $
b_i + a_{\tau ^{-1} (i)} +b _{\mu^{-1}(i)} \equiv a_i$ for $i=
(k-1)r+1, \cdots, (k-1)r+r.$ Let $i_1 =: (k-1)r+1, $ $i_2 =:
(k-1)r+2, $, $\cdots,$ $i_r =: (k-1)r+r, $ and $i_j' =: \tau (i_j)$
for $1\le j \le r.$ It is clear that (\ref {e6.1.2}) holds.
Consequently, it follows from the proof of part (a) that $(b, \tau )\in W(D_n)$.

(d). If the sign cycle $(b, \tau ) \in W(B_n)^{(a, \mu)}$ with $\tau =
B_{h, r}$, then $b(\tau \cdot a) (\mu \cdot b) \equiv a$. This
implies that $ b_i + a_{\tau ^{-1} (i)} +b _{\mu^{-1}(i)} \equiv a_i$ for
$i= (h-1)r+1, \cdots, (h-1)r+r,$ which is equivalent to $b _{(h-1)r
+j} + a_{hr+j} + b_{(h-1)r+ j-1} \equiv a_{(h-1)r +j}$ and $b _{hr
+j} + a_{(h-1)r+j} + b_{hr+ j-1} \equiv a_{hr +j}$. By addition of
the two equations we have  $ (b _{(h-1)r +j} + b _{hr +j}) +
(a_{hr+j} + a_{(h-1)r+j}) + (b_{(h-1)r+ j-1} + b_{hr+ j-1}) \equiv
(a_{(h-1)r +j} + a_{hr +j})$, i.e. $(b _{(h-1)r +j} + b _{hr +j}) +
(b _{(h-1)r +j-1} + b _{hr +j-1})\equiv 0$. Consequently, signs of
the cycles of $((h-1)r+j, hr+j)$ and $((h-1)r+j-1, hr+j-1)$ in $(b,
\tau)$ are the same  for $j = 2, 3,\cdots r$. Furthermore, the signs
of all cycles of $((h-1)r+j, hr+j)$ in $(b, \tau)$  for $j =1, 2,
3,\cdots r$ are the  same. Considering that the number of cycles in
$\tau$ is even, we have $(b, \tau ) \in W(D_n).$
$\Box$

\begin {Lemma}\label {6.3} Let $ a, c \in A$ and $\tau , \mu \in \mathbb
S_n$. Then $a$ and $ c$ are conjugate in
$W(B_n)$ if and only if $a$ and $ c$ are conjugate in $W(D_n)$;
$\tau$ and $\mu$ are conjugate in $W(B_n)$ if and only if $\tau$ and
$\mu$ are conjugate in $W(D_n)$.

\end {Lemma}
\noindent {\bf Proof.} The second claim  is clear. If there exists $(b, \xi)
\in W(B_n) = A\rtimes \mathbb S_n$ such that $(b, \xi) (a) (b,
\xi)^{-1} = c$. Thus $\xi \cdot a = c$ and $\xi a \xi^{-1} = c$,
which implies that $a$ and $ c$ are conjugate in $W(D_n)$. $\Box$

\section*{Acknowledgments}
  Authors  were financially
supported by the Australian Research Council. S. Zhang thanks the
School of Mathematics and Physics, The University of Queensland, and
Department of Mathematics, the Hong Kong University of Science and
Technology, for hospitality. We sincerely thank the referee for his or her kind help.


\begin {thebibliography} {200}

\bibitem {AS98} N. Andruskiewitsch and H. J. Schneider,
Lifting of quantum linear spaces and pointed Hopf algebras of order
$p^3$,  {\em J. Alg.} {\bf 209} (1998) 645-691.

\bibitem {AS02} N. Andruskiewitsch and H. J. Schneider,
Pointed Hope algebras, in {\em New directions in Hopf algebras}, Math. Sci. Res. Inst.
Publ. {\bf 43} (Cambridge University Press, 2002), pp. 1-68. 


\bibitem {Gr00} M. Gra\~na, On Nichols algebras of low dimension,
 {\em Contemp. Math.}  {\bf 267}  (2000) 111-134.

\bibitem {AFGV} N. Andruskiewitsch,  F. Fantino, M. Gra\~na and
L.Vendramin, Finite-dimensional  pointed Hopf algebras with
alternating groups are trivial, {\em Ann. Mat. Pura Appl.} {\bf 190} (2011) 225-245.

\bibitem {AZ07} N. Andruskiewitsch and S. Zhang, On pointed Hopf
algebras associated to some conjugacy classes in $\mathbb S_n$,
{\em Proc. Amer. Math. Soc.} {\bf 135} (2007) 2723-2731.

\bibitem {AFZ08}  N. Andruskiewitsch, F. Fantino, S. Zhang,   On pointed
Hopf algebras associated with the symmetric groups,   {\em Manuscripta
Math.}  {\bf 128} (2009) 359-371.

\bibitem {Bo68} N. Bourbaki, {\em Groupes et alg\`ebres de Lie} (Hermann, Paris, 1968).

\bibitem {Sa01} B. E. Sagan, The symmetric group: representations, combinatorial
algorithms, and symmetric functions, {\em Graduate Texts in Mathematics} (Springer-Verlag, 2001).

\bibitem {He06} I. Heckenberger, The Weyl groupoid of a Nichols algebra
of diagonal type, {\em Invent. Math.} {\bf 164} (2006) 175-188.

\bibitem {Se77} J.-P. Serre, {\em Linear representations of finite groups} 
(Springer-Verlag, New York, 1977).

\bibitem {HS10} I. Heckenberger and H.-J. Schneider, 
Root systems and Weyl groupoids for  Nichols algebras,  {\em Proc. London Math.} {\bf 101} (2010) 623-654.

\bibitem {AF07} N. Andruskiewitsch and F. Fantino, On   pointed Hopf
algebras with alternating and dihedral groups, {\em Rev. Uni¨®n Mat.
Argent.}  {\bf 48}  (2007) 57-71. 



\bibitem  {Su78} M. Suzuki,  {\em Group theory I} (Springer-Verlag, New York, 1978).

\bibitem {Ca72}  R. W. Carter, Conjugacy classes in the Weyl
group, {\em Compositio Mathematica}, {\bf 25} (1972) 1-59.







\end {thebibliography}
\end {document}